\documentclass[12pt]{article}
\usepackage[margin=25mm,centering]{geometry}

\usepackage{lineno}
\newtheorem{remark}{{Remark}}[section]

\title{Fast optimization of viscosities for frequency-weighted damping of second-order systems}
\date{\small \textbf{Submitted:} April 8, 2021}

\author{
Nevena Jakov\v{c}evi\'{c} Stor\thanks{
Faculty of Electrical Engineering, Mechanical Engineering and
Naval Architecture, University of Split, Rudjera Bo\v{s}kovi\'{c}a
32, 21000 Split, Croatia,
\texttt{nevena@fesb.hr}}
\and
Tim Mitchell\thanks{Max Planck Institute for Dynamics of Complex Technical Systems,
Sandtorstr. 1, 39106 Magdeburg, Germany,
\texttt{mitchell@mpi-magdeburg.mpg.de}}
\and
Zoran Tomljanovi\'c\thanks{Department of Mathematics,
University Josip Juraj Strossmayer, Trg Ljudevita Gaja 6, 31000 Osijek, Croatia,
\texttt{ztomljan@mathos.hr}, \texttt{mugrica@mathos.hr}}
\and
Matea Ugrica\footnotemark[3]
 }

\usepackage{amsmath}
\usepackage{amsfonts}
\usepackage{amssymb}
\usepackage{latexsym}
\usepackage{graphicx}
\usepackage{epstopdf}
\usepackage[shortlabels]{enumitem}
\usepackage{pgfplots}
\usepackage{tikz}
\usepackage{xfrac}
\usepackage{footnote}
\usepackage{hyperref}
\usepackage{subfig}
\usepackage[capitalize]{cleveref}
\usepackage{mathtools}
\usepackage{xspace}
\usepackage{color}
\usepackage{verbatim}

\usetikzlibrary{calc,patterns,decorations.pathmorphing,decorations.markings}

\numberwithin{equation}{section}

\newcommand{\MATLAB}{MATLAB}
\DeclareMathOperator{\diag}{diag}

\newcommand{\prtl}[2]{\frac{\partial#1}{\partial#2}}

\def\bigO{\mathcal{O}}
\def\imagunit{\mathbf{i}}
\def\transsym{\mathsf{T}}
\newcommand{\tp}[1][]{^{{#1}\transsym}}

\def\R{\mathbb{R}}
\def\C{\mathbb{C}}

\def\nnR{\R_+}

\renewcommand{\Re}{\mathrm{Re}\,}
\renewcommand{\Im}{\mathrm{Im}\,}
\DeclareMathOperator{\sgn}{sgn}

\usepackage{algorithm}
\usepackage{algpseudocode}
\usepackage{etoolbox}

\makeatletter
\AfterEndEnvironment{algorithm}{\let\@algcomment\relax}
\AtEndEnvironment{algorithm}{\kern2pt\hrule\relax\vskip3pt\@algcomment}
\let\@algcomment\relax
\newcommand\algcomment[1]{\def\@algcomment{\footnotesize \sc{Note: \it#1}}}

\renewcommand\fs@ruled{\def\@fs@cfont{\bfseries}\let\@fs@capt\floatc@ruled
  \def\@fs@pre{\hrule height.8pt depth0pt \kern2pt}%
  \def\@fs@post{}%
  \def\@fs@mid{\kern2pt\hrule\kern2pt}%
  \let\@fs@iftopcapt\iftrue}
\makeatother

\def\figsizes{8.0cm}	

\begin{document}

\maketitle

\begin{abstract}

We consider frequency-weighted damping optimization for vibrating systems described by a second-order differential equation.
The goal is to determine viscosity values such that eigenvalues are kept away from certain undesirable areas on the imaginary axis.
To this end, we present two complementary techniques.
First, we propose new frameworks using nonsmooth constrained optimization problems,
whose solutions both damp undesirable frequency bands and maintain stability of the system.
 These frameworks also allow  us to weight which frequency bands are the most important to damp.
Second, we also propose a fast new eigensolver for the structured quadratic eigenvalue problems that appear in such vibrating systems.
In order to be efficient, our new eigensolver exploits special properties of diagonal-plus-rank-one complex symmetric matrices,
which we leverage by showing how each quadratic eigenvalue problem can be transformed
into a short sequence of such linear eigenvalue problems.  The result is an eigensolver that is
substantially faster than standard techniques.
By combining this new solver with our new optimization frameworks,
we obtain our overall algorithm for fast computation of optimal viscosities.
The efficiency and performance of our new methods are verified and illustrated on several numerical examples.
\end{abstract}


\section{Introduction}
\label{intro}
Consider a vibrational mechanical system described by the second-order differential equation
\begin{equation}
	\label{sys_eq}
	M\ddot{q}(t)+C(v)\dot{q}(t)+Kq(t)=0,
\end{equation}
where $M, C(v), K \in \mathbb{R}^{n\times n}$ are all symmetric positive definite matrices, respectively representing mass, damping, and stiffness, and the damping matrix $C(v)$ depends on $r$ nonnegative viscosity parameters, i.e., $v \in \nnR^r$,
where $\nnR^r$ is the set of $r$-dimensional vectors with real nonnegative entries.
We assume that the number of damping parameters is small, i.e., $r \ll n$, as is typical in practice,
and that $C(v)$ has the following form:
\begin{equation}
	\label{damp_matr}
	C(v)=C_\mathrm{int} + G\,\diag(v_1,\dots,v_r)
	G\tp = C_\mathrm{int}+\sum\limits_{j=1}^{r}v_j g_jg_j\tp,
\end{equation}
where $C_\mathrm{int}$ represents internal damping, $G \in \R^{n \times r}$  describes the geometry of damping positions,
and $g_j$ denotes the $j$th column of $G$.
Internal damping can be modeled in various ways, e.g., Rayleigh (or classical) damping, where $C_\mathrm{int}=\alpha  M+\beta K$, where $\alpha, \beta \geq 0$.
In this paper, we focus another convention that is widely used, namely that the internal damping is a small multiple of the critical damping,~i.e.,
\begin{equation}
	\label{damp_int}
	C_\mathrm{int}=\alpha M^{\frac 1 2}\sqrt{ M^{-\frac 1 2}K  M^{-\frac 1 2}} M^{\frac 1 2},
\end{equation}
where $\alpha > 0$.  In this case, $C_\mathrm{int}$ is symmetric positive definite.
For more details on critical damping, see \cite{Ves11,KuzTT12,Adh06}.

The second-order differential equation \eqref{sys_eq} of course is associated with the quadratic eigenvalue problem (QEP)
\begin{equation}
	\label{eq:qep_eq}
	(\lambda(v)^2M+\lambda(v) C(v)+K)x(v)=0.
\end{equation}
For a given vector $v$ specifying the viscosity parameters,
let $\Lambda(v)$ denote the spectrum of~\eqref{eq:qep_eq}.
Each eigenvalue $\lambda(v)$ corresponds to a natural frequency of the system \eqref{sys_eq}, i.e., a frequency on which the system prefers to vibrate. Vibrations can be increased if the system is excited by an external force whose frequencies are close to the natural frequencies.  All frequencies that can significantly excite the system are called undesirable frequencies.

One approach to minimizing the impact of external forces is  damping optimization, which has been widely studied in the last few decades. In the most general context, given mass and stiffness matrices, the problem is to determine a damping matrix such that unwanted vibrations decay as fast as possible. This requires specifying an objective function to be optimized, and the  choice of objective function strongly depends on the application and desired outcome.   An overview of different damping optimization criteria can be found in~\cite{Ves11}. For the non-homogeneous case, where the system is additionally excited, damping optimization has also been studied in \cite{KuzTT16,Ves11}. In \cite{KuzTT16}, the authors consider energy over arbitrary time, while \cite{Ves11} considers the case where the excitation function is periodic.
For multiple input, multiple output systems, one can also optimize damping in second-order systems by minimizing standard systems norms, such as the $\mathcal{H}_2$ or $\mathcal{H}_\infty$ norms; see \cite{morBeaGT20,morBenKTetal16,BlaCGetal12,morTomBG18,morTomV20}.

In our setting, \eqref{sys_eq} corresponds to the homogeneous case for which one can consider optimizing the total average energy in various ways; see \cite{morBenTT11a,morBenTT13,CoxNRetal04,Ves11}.
One can also use eigenvalue-based criteria to damp resonant frequencies, i.e., by directly altering the spectrum of \eqref{eq:qep_eq},
as has been considered in \cite{morGraMQetal16} where the spectral abscissa criterion is minimized.
The spectral abscissa of \eqref{eq:qep_eq} is defined
\begin{equation}
  	\alpha_{\mathrm{MCK}}(v) =\max\limits_{\lambda(v) \in \Lambda(v)} \Re \lambda(v).
\end{equation}
Given some $v \in \R^r$, the system \eqref{sys_eq} is asymptotically stable  if and only if all the eigenvalues of the corresponding eigenvalue problem \eqref{eq:qep_eq} are in the open left half-plane, i.e., \mbox{$\alpha_\mathrm{MCK}(v) <0$}.
Note that under our assumption that $M, C(v), K$ are all symmetric positive definite matrices, and where $v \in \nnR^r$,
the system \eqref{sys_eq} is asymptotically stable; for more details, see \cite{TisM01}.

We consider the frequency isolation problem where viscosities are optimized in order to keep eigenvalues away from the certain undesirable areas on the imaginary axis, i.e., away from undesirable resonant frequency bands that are known \emph{a priori}.
This variation of the frequency isolation problem has been studied in several works.  In \cite{Jos92},
a Newton-type method for structures vibrating at low frequencies was proposed, while
a less costly inverse eigenvalue method was presented in \cite{EgaKS02},
where a target spectrum, which avoids an undesirable resonance band, is fixed in advance.
Meanwhile, \cite{MorE16} considered the frequency isolation problem for undamped vibrational systems where there is no $C(v)$ in equation~\eqref{eq:qep_eq}, or equivalently, $C(v)$ is always zero.
In damping optimization where $C(v)$ is present, avoiding undesirable frequency bands
can be achieved by either choosing damping positions (by optimizing matrix $G$)
or by damping viscosities (by optimizing $v \in \nnR^r$) or doing both simultaneously.
Computing the optimal damping positions is a very challenging problem and there is no efficient algorithm for it,
though some heuristics can be found in, e.g., \cite{KanPTetal19}.
One approach to determining optimal damping positions is ``direct" brute force,
where all possible damping configurations are considered and viscosities are optimized for each configuration.
In any case, while optimization of damping positions is a challenging and a very important question in and of itself,
in this paper we focus on accelerating this overall process via proposing faster methods for
viscosity optimization for fixed damping positions.
Therefore, in the our algorithms here, we consider that the matrix $G$ specifying the damping positions is fixed,
but we have in mind that, in practice, viscosity optimization would be applied over many different configurations of damping positions.

In this paper, we aim to accelerate such damping-based approaches for frequency isolation
via proposing new fast techniques for the important subproblem of determining optimal damping viscosities for a given configuration of damping positions.
More specifically, given a general system \eqref{sys_eq}, where $C(v)$ has the form given in \eqref{damp_matr}, the internal damping matrix $C_\mathrm{int}$ is given by \eqref{damp_int}, and the matrix $G$ specifying the damping positions is fixed, we consider the problem of optimizing the viscosity parameters $v \in \nnR^r$ so that
the eigenvalues of \eqref{eq:qep_eq} are kept away from undesirable resonant bands as much as possible.
Our contribution here actually consists of two complementary new techniques.
First, we propose two related nonsmooth but continuous constrained optimization problems as new models
for the frequency isolation problem and show how solutions can be computed via gradient information
and recent advances in nonsmooth constrained optimization.
When our new problems are solved, their solutions
provide viscosity parameters which maintain stability of the system and damp
user-defined undesirable frequency bands.
In addition to specifying the number of frequency bands and their respective widths, the undesirable bands
can also be weighted in order to prioritize which are most critical to damp.
Second, as the cost of our optimization process is actually dominated by solving a sequence of related QEPs,
where $C(v)$ is changing as the viscosity parameters are optimized,
we also propose a fast algorithm to solve this sequence of QEPs.
Our method, which is many times faster than using standard eigensolvers for QEPs
and can be considered an extension of \cite{JakST20} for computing eigenvalues
of  diagonal-plus-rank-one (DPR1) complex symmetric (DPR1Csym) matrices,
works by exploiting the fact that changing the viscosity parameters is equivalent to making a
low-rank update to $C(v)$.
Since such structure is not inherent to our problem, we expect that
our technique for solving such sequences of QEPs could be quite beneficial in other applications as well.

The paper is organized as follows. In \cref{frame_formulation}, we motivate and establish our two new models for
the frequency isolation problem, explain their properties, and discuss how to compute solutions of them.
Then, in \cref{QEPs_dpr1}, we present our new approach for efficiently
solving the corresponding sequence of QEPs that arises during optimization (using either of our new models discussed in
the previous section).  We show how both eigenvalues and eigenvectors of the QEPs can be computed, as each
are needed in our optimization-based approach.
In \cref{sec:damp_opt}, we present our full algorithm for damping optimization by combining our aforementioned components from \cref{frame_formulation,QEPs_dpr1}.
Finally, we validate our new techniques and compare our two models for frequency-weighted damping in \cref{nume_ex}.

\section{New frameworks for frequency-weighted damping} \label{frame_formulation}
Consider how the eigenvalues $\Lambda(v)$ of \eqref{eq:qep_eq} evolve
as the viscosities parameters $v$ are changed, and suppose
$\omega \in \R$ is an undesirable frequency, i.e., we wish to keep the spectrum of \eqref{eq:qep_eq}
away from $\imagunit \omega$ on the imaginary axis.
Since eigenvalues with imaginary parts close to $\omega$ can also be undesirable,
we thus consider the undesirable frequency band $[\omega - b,\omega +b]$ about $\omega$
for some given $b > 0$.
In order to minimize the impact of eigenvalues of \eqref{eq:qep_eq} in this frequency band,
it is tempting to consider solving the optimization problem:
\begin{equation}
	\label{eq:naive_opt}
	\begin{aligned}
	\min_{v \in \R^r} \quad & \max \{ \Re \lambda(v) :  \lambda(v) \in \Lambda(v) \text{ and } \Im \lambda \in [\omega - b,\omega +b] \} \\
	\text{s.t.} \quad & \alpha_\mathrm{MCK}(v) \leq \texttt{tol}_\mathrm{sa} \text{ for some } \texttt{tol}_\mathrm{sa} < 0, \\
				& v_j \geq 0 \text{ for } j = 1,\ldots,r,
	\end{aligned}
\end{equation}
which would act to push all the eigenvalues of $\Lambda(v)$ with imaginary parts in $[\omega - b,\omega +b]$
as far to the left as possible while still maintaining asymptotic stability of the system
and physically realistic, i.e., nonnegative, viscosities.
Alternatively, one might consider swapping the objective function and the stability constraint in \eqref{eq:naive_opt},
i.e., minimize the spectral abscissa as much as possible subject to the constraint that
the eigenvalues in the frequency band  $[\omega - b,\omega +b]$ are all kept at least
some fixed distance to the left of the imaginary axis (and again enforcing nonnegative viscosities).
However, these two related optimization problems are rather difficult to solve
as the function being minimized in \eqref{eq:naive_opt} is actually discontinuous.
In general, this function has jump discontinuities
whenever a rightmost eigenvalue that attains the maximum
leaves the horizontal strip in the complex plane defined by $\omega$ and $b$, or vice versa,
when a new eigenvalue enters this region to become a rightmost eigenvalue in this strip,
and these discontinuities are typically not uncommon.

To overcome this problem, in this section we propose two alternatives to \eqref{eq:naive_opt}
where continuity is maintained and so our new optimization problems for frequency isolation are much more practical to solve.
This allows us to use continuous optimization techniques to compute viscosity values
such that eigenvalues are kept away from an undesirable frequency band defined by $\omega$ and $b$.
In fact, as we will soon explain, our distance function can be used for different undesirable frequency bands simultaneously.

\subsection{Approach 1}
\label{sec:approach1}
Let the tuple $E = (a,b,c)$ denote the axis-aligned ellipse
\begin{equation}
	\label{ellipse1}
	 \frac{(x - \Re c)^2}{a^2}+\frac{(y- \Im c)^2}{b^2}=1,
\end{equation}
where $a,b > 0$ respectively denote the semi-major and -minor axes and $c \in \C$ is the center of the ellipse.
Identifying $\R^2$ with $\C$,
consider the following algebraic distance $d : \C \mapsto [0,\infty)$ of a point $z \in C$ to this ellipse, i.e.,
\begin{equation}
	\label{eq:dist_fun}
	d(z;E) \coloneqq
	\frac{(\Re (z - c))^2}{a^2} + \frac{( \Im (z - c))^2}{b^2}.
\end{equation}
Thus, $d(z;E) > 1$ when $z$ is outside of the ellipse, $d(z;E) \in [0,1)$ when $z$ is inside the ellipse,
and $d(z;E)=1$ when $z$ is on the ellipse, i.e., $z = x +  \imagunit y$ satisfies \eqref{ellipse1}.

Now suppose that $\omega \geq 0$ and $b > 0$ specify an undesirable frequency band $[\omega - b,\omega + b]$.
Given some $a > 0$, we can measure the distance to a point in the complex plane to the interval $\imagunit[\omega - b,\omega + b]$
on the imaginary axis via $d(z;E)$ for $E = (a,b,\imagunit \omega)$, i.e., the ellipse \eqref{ellipse1} centered at $\imagunit \omega$ on the imaginary axis.
If $z$ is such that $\Im z \in (\omega - b,\omega + b)$,
then the larger we make the value of $a$, the further $z$ must be to the left or right of the minor axis of the ellipse given by $E$
in order for $d(z;E) > 1$ to hold.
Thus, as a continuous measure of the distance of the spectrum of \eqref{eq:qep_eq} to the undesirable frequency
$[\omega - b,\omega + b]$, we define
\begin{equation}
	\label{eq:dist_fun_single}
	d_{\Lambda,E}(v) \coloneqq \min \{ d(\lambda(v); E) : \lambda(v) \in \Lambda(v) \},
\end{equation}
where $E = (a,b,\imagunit \omega)$.
Function  $d_{\Lambda,E}(v) > 1$ when all the eigenvalues of $\Lambda(v)$ are outside the given ellipse,
$d_{\Lambda,E}(v) \in [0,1)$ when one or more eigenvalues are inside this ellipse,
and $d_{\Lambda,E}(v) = 1$ when at least one eigenvalue is on this ellipse and none are inside.
The specific value of the semi-major axis $a$ determines the importance of the undesirable frequency band
by dictating how far away eigenvalues should be from the interval $\imagunit[\omega - b,\omega + b]$,
where eigenvalues with imaginary parts closer to $\omega$ are weighted more, i.e., must be further away.
When multiple undesirable frequency bands are specified by
frequencies $\{\omega_1,\ldots,\omega_k\}$ and associated (half) bandwidths $\{b_1,\ldots,b_k\}$,
their relative importance can be determined by providing different semi-major axis values $\{a_1,\ldots,a_k\}$,
with $\omega_j \geq 0$ and $a_j,b_j > 0$ for all \mbox{$j=1,\ldots,k$}.
Thus, we generalize \eqref{eq:dist_fun_single} to measuring the distance of the spectrum $\Lambda(v)$
to the $k$ undesirable frequency bands by defining
\begin{equation}
	\label{eq:dist_fun_mult}
	d_{\Lambda,\mathcal{E}}(v) \coloneqq \min \{ d_{\Lambda,E_j}(v) : E_j \in \mathcal{E}  \},
\end{equation}
where $E_j \coloneqq (a_j,b_j,\imagunit \omega_j)$ is the tuple defining the $j$th axis-aligned ellipse for the $j$th
undesirable frequency band $[\omega_j - b_j, \omega_j + b_j]$ with relative importance $a_j > 0$
and $\mathcal{E} \coloneqq \{E_1,\ldots,E_k\}$ is the set of $k$ corresponding ellipses.

Using \eqref{eq:dist_fun_mult}, we now present our first new model for the frequency isolation problem:
\begin{equation}
	\label{eq:model1}
	\hspace*{-2cm}
	\begin{aligned}
	\textbf{Model 1:} \qquad
	\min_{v \in \R^r} \quad & \alpha_\mathrm{MCK}(v) \\
	\text{s.t.} \quad & d_{\Lambda,\mathcal{E}}(v) \geq 1, \\	
				& \alpha_\mathrm{MCK}(v) \leq \texttt{tol}_\mathrm{sa} \ \text{ for some } \ \texttt{tol}_\mathrm{sa} < 0, \\
				& v_j \geq 0 \text{ for } j = 1,\ldots,r,
	\end{aligned}
\end{equation}
i.e., minimize the spectral abscissa as much as possible subject to the respective constraints
that all the eigenvalues of \eqref{eq:qep_eq} are outside of the ellipses defined by $\mathcal{E}$,
the system is asymptotically stable, and the viscosities are nonnegative.
Although the spectral abscissa is being minimized in \eqref{eq:model1},
note that the additional constraint that the spectral abscissa be negative
is necessary.  There are multiple reasons for this.
First, not all optimization solvers iterate only over the feasible set, and so negative viscosities
may be encountered during optimization, which in turn may make $\alpha_\mathrm{MCK}(v)$ positive.
Second, satisfying $d_{\Lambda,\mathcal{E}}(v) \geq 1$ is not equivalent
to satisfying stability, as $d_{\Lambda,\mathcal{E}}(v) \geq 1$ can hold even if
all the eigenvalues were to be in the right half-plane.
Third, $\alpha_\mathrm{MCK}(v)$ may have stationary points where $\alpha_\mathrm{MCK}(v) \geq 0$ holds,
and so a feasible minimizer of \eqref{eq:model1} without this stability constraint
would not necessarily result in an asymptotically stable system.

While $\alpha_\mathrm{MCK}(v)$ and $d_{\Lambda,\mathcal{E}}(v)$ in \eqref{eq:model1} are nonsmooth functions,
they are at least continuous (unlike the objective function in \eqref{eq:naive_opt}).
As there has been significant progress recently in developing effective solvers for nonsmooth constrained optimization,
e.g., \cite{CurO12,CurMO17},
where the functions are continuous but their nonsmoothness is restricted to a set of measure zero, as is typical,
it is reasonable to apply such techniques in order to compute minimizers of \eqref{eq:model1}.
We describe the details of how this is done later on and
for now make some additional general comments about \eqref{eq:model1}.
Since  $\alpha_\mathrm{MCK}(v)$ and $d_{\Lambda,\mathcal{E}}(v)$ will typically be nonconvex
and \eqref{eq:model1} may have infeasible stationary points,
we cannot necessarily expect to find a globally optimally solution to \eqref{eq:model1},
and solvers may also sometimes converge to infeasible points.
However, in applications, local optimally solutions are often sufficient and
provide meaningful improvements in performance over non-optimized configurations.
Moreover, both of these problems can typically be mitigated merely by computing multiple solutions to \eqref{eq:model1}
via initializing a solver from many different starting points and taking the best of the resulting solutions.
Note that the choice of $\mathcal{E}$ depends on the application and is fixed before optimization commences.
However, if the $a_j$ values are chosen too aggressively (too large), there is no guarantee that \eqref{eq:model1}
will have any feasible solutions.
Thus, we now propose a second new model as an alternative and which avoids this issue.

\subsection{Approach 2}
Again consider a single axis-aligned ellipse \eqref{ellipse1} defined by tuple $E = (a,b,c)$, where
$a,b > 0$ and $c = \eta + \imagunit \omega$ with $\eta, \omega \geq 0$ (so $c$ is not necessarily on the imaginary axis),
and suppose that $b$ and $c$ are fixed but $a$ may be varied.
Then, given a point $z \in \C$,
consider the largest we can make this ellipse, by changing the length of its major axis, such that
$z$ is not inside the ellipse.  For $z = x + \imagunit y$ and $c=\eta + \imagunit \omega$,
solving \eqref{ellipse1} for yields that this largest possible value for the semi-major axis is:
\begin{equation}
	\label{eq:major_fun}
	a(z;E) \coloneqq
	\begin{cases}
		\frac{ b | \Re z - \eta |}{\sqrt{ b^2 - (\Im z - \omega)^2}}, & \quad \text{if } \Im z \in (\omega - b, \omega + b), \\
		\infty & \quad \text{otherwise},
	\end{cases}
\end{equation}
where the largest possible semi-major axis value is infinite when the point $z$ is not directly to the left or right of the ellipse,
i.e., $\Im z \not\in (\omega - b, \omega + b)$.  Note that this convention is consistent even when
$z = \eta + \imagunit(\omega \pm b)$, i.e., one of the endpoints of the minor axis, since in this case, $z$ can never be inside the
ellipse no matter how large the major axis is.
While $a(z;E)$ is determined only by $b$ and $c$ from the tuple $E$,
we continue to use the tuple $E=(a,b,c)$ for notational and conceptual consistency with \cref{sec:approach1},
but when the value of $a$ is not fixed, we will instead write $E = (\sim,b,c)$.

As a function of $z$, $a(z;E)$ is real valued and always nonnegative.
Note that $a(z;E)$ is continuous wherever $\Re z \neq \eta$,
since then the ratio in \eqref{eq:major_fun} continuously goes to infinity as \mbox{$\Im z \in (\omega -b, \omega + b)$}
approaches $\omega \pm b$.
When $\Re z = \eta$ holds,
$a(z;E)$ only has two discontinuities, as in this case, the numerator in $a(z;E)$ is zero,
and so $a(z;E)$ has a jump between zero and infinity at $z = \eta + \imagunit (\omega \pm b)$.
Relative to \eqref{eq:naive_opt}, where the discontinuities can be common and negatively impact solvers,
the only two discontinuities of $a(z;E)$ are relatively benign as they typically will not be encountered.
Moreover, by a modification which we will explain momentarily,
we can in fact completely remove these discontinuities from the optimization problem.

Now considering the spectrum $\Lambda(v)$ of \eqref{eq:qep_eq},
we can use $a(z;E)$ to determine how the largest value of the semi-major axis of the ellipse $E=(\sim,b,\eta + \imagunit \omega)$,
such that none of the eigenvalues are inside it, varies with respect to the viscosities $v$ changing:
\begin{equation}
	a_{\Lambda,E}(v) \coloneqq
	\min \{ a(\lambda(v); E) : \lambda(v) \in \Lambda(v) \}.
\end{equation}
Clearly, function $a_{\Lambda,E}(v)$ inherits the properties of $a(v;E)$ discussed above, but
since $\eta \geq 0$, note that $a_{\Lambda,E}(v)$ is continuous at $\hat v$ if the system is asymptotically stable
for $\hat v$.  Furthermore, $a_{\Lambda,E}(v)$ is smooth at point $\hat v$
if there is only a single eigenvalue (excluding conjugacy) on the ellipse given by
$E = (a_{\Lambda,E}(\hat v), b, \eta + \imagunit \omega)$ and this eigenvalue is simple.
Thus, to damp the frequency band $[\omega - b, \omega +b]$, we could consider solving,
\begin{equation}
	\label{eq:model2_disc}
	\begin{aligned}
	\max_{v \in \R^r} \quad & a_{\Lambda,E}(v) \\
	\text{s.t.} \quad & \alpha_\mathrm{MCK}(v) \leq \texttt{tol}_\mathrm{sa} \ \text{ for some } \ \texttt{tol}_\mathrm{sa} < 0, \\
				& v_j \geq 0 \text{ for } j = 1,\ldots,r.
	\end{aligned}
\end{equation}
where maximizing $a_{\Lambda,E}(v)$ acts to push all eigenvalues directly
to the left of the interval \mbox{$\imagunit[\omega - b, \omega +b]$} as far to the left as possible.
While \eqref{eq:model2_disc} is still discontinuous,
encountering the two discontinuities of $a_{\Lambda,E}(v)$ during optimization is unlikely;
not only do the discontinuities occur off of the feasible set, i.e., when the system is not stable,
they require that an eigenvalue passes through  $\eta + \imagunit (\omega \pm b)$ \emph{exactly}
in order to occur.
On the other hand, many optimization solvers do explore the infeasible set during optimization,
and even though function $a_{\Lambda,E}(v)$ is otherwise continuous,
it can nevertheless have arbitrarily high growth when there exists an eigenvalue $\lambda(v)$ with $\Re \lambda(v) \approx \eta$
and $\Im \lambda(v) \in (\omega - b, \omega + b)$ approaches the endpoints of this interval.
As such, before extending \eqref{eq:model2_disc} to the case of multiple ellipses, i.e., multiple frequency bands to damp,
we first propose modifying \eqref{eq:model2_disc} via a barrier function.

The core idea of introducing a barrier function is to alter \eqref{eq:model2_disc}
such that viscosities which cause $\alpha_\mathrm{MCK}(v)$
to get close to $\eta$ will be increasingly penalized, to the point that
optimization will never allow a configuration $\hat v$  to
be accepted as an iterate where $\alpha_\mathrm{MCK}(\hat v) \geq \eta$ holds.
We do this by modifying the objective function such that it goes to negative infinity as
$\alpha_\mathrm{MCK}(v)$ goes from $\texttt{tol}_\mathrm{sa}$ to $\eta$.
Since optimization can never accept points where the objective function is infinite,
this barrier guarantees that
points where $a_{\Lambda,E}(v)$ is discontinuous are never encountered.
Furthermore, accepting points where $a_{\Lambda,E}(v)$ is nearly
discontinuous will also be heavily discouraged, as the barrier-modified objective function that we are trying to maximize
quickly goes to negative infinity as $\alpha_\mathrm{MCK}(v)$ increases beyond $\texttt{tol}_\mathrm{sa}$.
However, such a barrier function should not introduce
any new discontinuities or nonsmooth points of its own, nor should
it alter the objective function where $\alpha_\mathrm{MCK}(v) \leq \mathrm{tol}_\mathrm{sa}$ holds,
as all of these things could make optimization more difficult.
We construct our barrier function out of a cubic polynomial
and log-based function that are specifically crafted to meet these goals.

Given real scalars $y_1 < y_2$ and a continuous function $f : \R^r \mapsto \R$,
we define the following generic barrier function
\begin{equation}
	\label{eq:barrier}
	\beta(f(x); y_1,y_2) \coloneqq
	\begin{cases}
	0 & \text{if } f(x) \leq y_1, \\
	\tau_1(f(x) - y_1)^3 + \tau_2(f(x) - y_1)^2 & \text{if } f(x) \in (y_1,y],  \\
	-\log \left( \frac{y_2 - f(x)}{y_2 - y} \right) + h & \text{if } f(x) \in (y,y_2),  \\
	\infty & \text{otherwise,}
	\end{cases}
\end{equation}
where $y \in (y_1,y_2)$, $h > 0$, and
\begin{equation}
	\label{eq:barrier_coeffs}
	\tau_1 \coloneqq \frac{(2h+1)y - y_1 - 2hy_2}{(y_2 - y)(y - y_1)^3}
	\quad \text{and} \quad
	\tau_2 \coloneqq \frac{y_1 + 3hy_2 - (3h+1)y}{(y_2 - y)(y - y_1)^2}.
\end{equation}
Thus, as $f(x)$ goes from $y_1$ to $y_2$, our barrier function $\beta(f(x); y_1,y_2)$
goes from zero to infinity.
The constants $\tau_1, \tau_2 \in \R$ are specifically
chosen so that the value of $\beta(f(x); y_1,y_2)$ always varies continuously
and $\nabla \beta(f(x); y_1,y_2)$ is continuous
wherever $\nabla f(x)$ is continuous, with $\| \nabla \beta(f(x); y_1,y_2) \| = 0$ if
$f(x) = y_1$.
The continuity of the gradients can be verified by differentiating the component functions in \eqref{eq:barrier},
which are shown later in \eqref{eq:barrier_grad}.
In other words, $\beta(f(x); y_1,y_2)$ realizes our goals stated above, as it is a continuous barrier function
that can be added to any objective function without introducing any nonsmooth points of its own,
i.e., points where the gradient is not defined.
However, if $\hat x$ is a nonsmooth point of $f(x)$ with $f(\hat x) \in (y_1,y_2)$,
then naturally $\beta(f(x); y_1,y_2)$ must also be nonsmooth at~$\hat x$.

The values $y$ and $h$ determine exactly where $\beta(f(x); y_1,y_2)$ switches between
its second and third cases, i.e., where the cubic polynomial meets the $\log$-based function when \mbox{$f(x) = y$} and
$\beta(f(x); y_1,y_2) =h$.  While the latter monotonically increases with respect
to $f(x)$ increasing, note that this is not necessarily guaranteed for the cubic polynomial.
However, this can be enforced with a judicious choice of $y$.
For example, if we choose to set $\tau_2=0$, and consider $f(x)=x$ (so $x \in \R$),
then the cubic polynomial and its first derivative are always increasing for all $y > y_1$.
We can then simply solve for $y$ by setting the numerator of $\tau_2$ in \eqref{eq:barrier_coeffs} equal to zero,
which yields $y = y_1 + \tfrac{3h}{3h + 1}\delta$, where $\delta = y_2 - y_1 > 0$.

We now modify and extend \eqref{eq:model2_disc} to respectively make it continuous via our barrier function
and support damping multiple frequency bands.
Suppose multiple undesirable frequency bands are specified by
frequencies $\{\omega_1,\ldots,\omega_k\}$ and associated (half) bandwidths $\{b_1,\ldots,b_k\}$,
with their relative importance determined by $\{\phi_1,\ldots,\phi_k\}$,
where $\omega_j \geq 0$, $b_j > 0$, and $\phi_j \in (0,1]$ for all \mbox{$j=1,\ldots,k$}.
Then given some $\eta \geq 0$ and the corresponding ellipses $\{E_1, \ldots, E_k\}$ with $E_j = (\sim,b_j,\eta + \imagunit \omega_j)$,
our second model for optimizing viscosities is
\begin{equation}
	\label{eq:model2}
	\hspace*{-1cm}
	\begin{aligned}
	\textbf{Model 2:} \qquad
	\max_{v \in \R^r} \quad
		& \left(\sum_{j = 1}^k \phi_j \cdot \min \{a_{\Lambda,E_j}(v), m_j \} \right)
			- \beta(\alpha_\mathrm{MCK}(v); \texttt{tol}_\mathrm{sa}, \eta) \\
	\text{s.t.} \quad & \alpha_\mathrm{MCK}(v) \leq \texttt{tol}_\mathrm{sa} \text{ for some } \texttt{tol}_\mathrm{sa} < 0, \\
				& v_j \geq 0 \text{ for } j = 1,\ldots,r,
	\end{aligned}
\end{equation}
where we use $h \coloneqq 1$ and $y \coloneqq \texttt{tol}_\mathrm{sa} + \tfrac{3}{4}(\eta - \texttt{tol}_\mathrm{sa})$ for
our barrier function \eqref{eq:barrier}
and $m_j > 0$ is a fixed scalar denoting a desired upper bound on the damping of the $j$th frequency band,
i.e., a point at which the band can be considered sufficiently damped.
The sum in the objective function of \eqref{eq:model2} acts to push all eigenvalues to the left
of the intervals $\imagunit [\omega_j - b_j, \omega_j + b_j]$
farther to the left, namely, by trying to increase each of the semi-major axis values of the ellipses
(while still having no eigenvalues inside them) as much as possible or until they are at least as large as the respective $m_j$ values.
The inclusion of the finite $m_j$ scalars prevent optimization terminating due to
one of the $a_{\Lambda,E_j}(v)$ functions becoming infinite, which happens
if all the eigenvalues can be moved completely outside of one or more of the specified frequency bands.
This can be undesirable because when this happens, the other frequency bands may or may not be well optimized.
Using $\min \{a_{\Lambda,E_j}(v), m_j \}$ prevents this from occurring, and so all the frequency bands
will continue to be optimized.
Meanwhile, the $\phi_j$ scalars allow one to balance which frequency bands should be given the most emphasis
during optimization.
By construction, our barrier function only has an effect when the spectral abscissa
stability constraint is violated, and so it does not modify our objective function on the feasible set.
As $\beta(\alpha_\mathrm{MCK}(v); \texttt{tol}_\mathrm{sa}, \eta)$
goes continuously to infinity as the spectral abscissa approaches $\eta$,
the discontinuities of $a_{\Lambda,E_j}(v)$ functions can never be encountered
and having eigenvalues with real parts close to $\eta$ is strongly penalized,
which helps to avoid regions where $a_{\Lambda,E_j}(v)$ is close to being discontinuous.
Compared to our first model in \cref{sec:approach1},
we have introduced the parameter $\eta \geq 0$ here so that if desired, the distance between
being stable to tolerance and the discontinuities of $a_{\Lambda,E_j}(v)$ can be increased
by shifting all the ellipses to the right.

\subsection{Solving our optimization problems}
\label{sec:opt_solve}
A key goal realized by our new constrained optimization problems for frequency-weighted damping,
respectively given in \eqref{eq:model1} and \eqref{eq:model2}, is that they are both continuous,
unlike the formulation we first considered in \eqref{eq:naive_opt}.
Consequently, as mentioned earlier,
we thus can consider computing solutions to  \eqref{eq:model1} and \eqref{eq:model2} using recent gradient-based solvers
for continuous nonsmooth constrained optimization, where the nonsmoothness of the functions are limited
to a set of zero.
This is appealing because such gradient-based nonsmooth solvers not only exhibit good performance in practice
but are also easy to use, as they only require that gradients be provided; see \cite[Section~6]{CurMO17} for some comparisons.
The necessary gradients exist because typically such methods only encounter the nonsmooth manifold in the limit,
and so while iterates may be arbitrarily close to nonsmooth points, they are nevertheless not nonsmooth points themselves.
Two possible gradient-based solvers for nonsmooth constrained optimization are the open-source software packages
SQP-GS \cite{CurO12} and
GRANSO: GRadient-based Algorithm for Non-Smooth Optimization \cite{CurMO17}.
For the purposes of this paper, we use GRANSO\footnote{Available at \url{https://gitlab.com/timmitchell/GRANSO/}.}
to compute solutions to \eqref{eq:model1} and \eqref{eq:model2},
partly because GRANSO is typically much faster than SQP-GS.
We now discuss how to compute the necessary gradients.

For our first approach, given by \eqref{eq:model1}, we need the gradient of the spectral abscissa
and $d_{\Lambda,\mathcal{E}}(v)$.  We begin with the former.
Let $\lambda(v)$ be an eigenvalue of \eqref{eq:qep_eq} with associated eigenvector $x(v)$.
Since $M$, $C(v)$, and $K$ are real symmetric matrices, $x(v)$ is also the left eigenvector for $\lambda(v)$.
Then given some $\hat v$, if $\lambda(\hat v)$ is a simple eigenvalue with eigenvector by $\hat x$,
by standard perturbation theory for eigenvalues we have that
\begin{equation}
\label{Lemma1:EigDeriv}
	\prtl{\lambda(v)}{v_j} \bigg|_{v = \hat v} =
	-\frac{\hat x^*\left( \lambda(\hat v) g_jg_j\tp \right) \hat x}{\hat x^*(2 \lambda(\hat v) M + C(\hat v)) \hat x}.
\end{equation}
Furthermore, if $\lambda(\hat v)$ is also an eigenvalue that attains the spectral abscissa
and there are no other eigenvalues with real part equal to $\Re \lambda(\hat v)$,
i.e., there are no ties (excluding conjugacy) for the spectral abscissa, then
\begin{equation}
	\label{eq:spec_abs_grad}
	\prtl{\alpha_\mathrm{MCK}(v)}{v_j} \bigg|_{v = \hat v} = - \Re \prtl{\lambda(v)}{v_j} \bigg|_{v = \hat v}.
\end{equation}
We now turn to $d_{\Lambda,\mathcal{E}}(v)$.
Given a single ellipse given by $E=(a,b,c)$,
consider $d(z(t); E)$ defined by \eqref{eq:dist_fun}, where $z(t)$ is a differentiable path with respect to the real scalar~$t$.
Then the derivative of $d(z(t); E)$ is
\begin{equation}
	\label{eq:dist_grad}
	d^{\, \prime}(z(t); E)
	= 2 \left( \frac{\Re (z(t) - c) \cdot \Re z^\prime(t)}{a^2}
	+\frac{\Im (z(t) - c) \cdot \Im z^\prime(t) }{b^2} \right).
\end{equation}
Now given $\hat v$, suppose there are no ties for the value of $d_{\Lambda,\mathcal{E}}(\hat v)$,
i.e., its value is attained by a single eigenvalue $\lambda(\hat v)$ and ellipse $E = (a,b,c) \in \mathcal{E}$,
with $\lambda(\hat v)$ being simple.
Then the gradient of $d_{\Lambda,\mathcal{E}}(v)$ at $\hat v$ exists,
and the partial derivative with respect to $v_j$ at $\hat v$ can be constructed via \eqref{eq:dist_grad}, where
$z(t)$ is replaced by $\lambda(\hat v)$ and $z^\prime(t)$ is replaced by the partial derivative of $\lambda(v)$
at $\hat v$
given in \eqref{Lemma1:EigDeriv}.

For our second approach, given by \eqref{eq:model2},
we have shown above how to obtain gradient of the spectral abscissa,
which leaves the objective function in \eqref{eq:model2}.
Given an ellipse defined by $E = (\sim,b,\eta + \imagunit \omega)$,
again consider $z(t)$ described above but additional suppose that $z(t) \in (\omega - b, \omega + b)$.
Then $a(z(t); E)$ cannot be infinite and its derivative is
\begin{equation}
	\label{eq:major_grad}
	a^\prime(z(t); E) =
	\frac{ b \sgn(\Re z(t) - \eta) \cdot \Re z^\prime(t) }{ (b^2 - (\Im z(t) - \omega)^2)^{\sfrac{1}{2}} }
	+
	\frac{ b | \Re z(t) - \eta| (\Im z(t) - \omega) \cdot \Im z^\prime(t) }{ (b^2 - (\Im z(t) - \omega)^2)^{\sfrac{3}{2}} }
\end{equation}
Now consider $a_{\Lambda,E_j}(v)$,
which is differentiable if there is a single eigenvalue (up to conjugacy) on the ellipse specified by $E_j$
and this eigenvalue is simple.
If these assumptions hold at $\hat v$ and this eigenvalue is $\lambda(\hat v)$,
then partial derivative with respect to $v_j$ of $a_{\Lambda,E_j}(v)$ at $\hat v$
is given by \eqref{eq:major_grad} with $z(t)$ and $z^\prime(t)$ are again replaced
using $\lambda(\hat v)$ and \eqref{Lemma1:EigDeriv}.
For the gradient of our barrier function, it suffices to show the derivative of
$\beta(f(x); y_1,y_2)$, where $f : \R \mapsto \R$ and is differentiable:
\begin{equation}
	\label{eq:barrier_grad}
	\beta^\prime(f(x); y_1,y_2) \coloneqq
	\begin{cases}
	0 & \text{if } f(x) \leq y_1, \\
	(3\tau_1(f(x) - y_1)^2 + 2\tau_2(f(x) - y_1))\cdot f^\prime(x) & \text{if } f(x) \in (y_1,y],  \\
	\frac{f^\prime(x)}{y_2 - f(x)} & \text{if } f(x) \in (y,y_2), \\
	\text{undefined} & \text{otherwise.}
	\end{cases}
\end{equation}
Then given $\hat v$ and assuming $\alpha_\mathrm{MCK}(v)$ is differentiable at $\hat v$,
the partial derivative with respect to $v_j$ of our barrier function in \eqref{eq:model2}
is given by \eqref{eq:barrier_grad}, where $f(x)$ and $f^\prime(x)$
are respectively replaced using $\alpha_\mathrm{MCK}(\hat v)$ and \eqref{eq:spec_abs_grad}, $y_1 = \texttt{tol}_\mathrm{sa}$,
and $y_2 = \eta$.

\section{Fast solution of QEPs with low-rank structure}
\label{QEPs_dpr1}
The most expensive part of our approaches proposed in \cref{frame_formulation} is successively computing the eigenvalues and eigenvectors of \eqref{eq:qep_eq} as the viscosities are optimized, i.e., as $C(v)$ is changed.
One possibility is to use \texttt{polyeig} in \MATLAB\ or \texttt{quadeig}; see \cite{HamMT13,TisM01}
for more details on these methods.
However, using either of these routines would mean that solving each QEP would require roughly the same amount of cubic work,
i.e., $\bigO(n^3)$, where we use the usual convention of treating eigenvalue computations as atomic operations.
In~\cite{Tas15}, Taslaman proposed a much faster eigensolver for QEPs~\eqref{eq:qep_eq},
where the damping matrix $C(v)$ is assumed to be low rank.
While the overall work complexity of Taslaman's algorithm is still cubic,
in experiments~\cite[section~5]{Tas15}, it was shown to be many times faster than \texttt{quadeig},
and its work can be separated into offline and online components, with the latter only doing $\bigO(n^2)$ work.
Taslaman's algorithm is based on Ehrlich-Aberth iterations, which  
rely on a good choice of a starting point for each eigenvalue and for which determination of stopping criteria
is often heuristic; for more details, see \cite{BinN13} and \cite{Tas15}.
Shortly thereafter, Taslaman's algorithm was extended by Benner and Deni{\ss}en~\cite{BenD15}
to systems where $C(v) = C_\mathrm{int} + C_\mathrm{ext}(v)$ may be full rank,
but critical damping is used for the internal damping matrix $C_\mathrm{int}$ and the 
external damping matrix $C_\mathrm{ext}(v)$ is low rank.

In this section, for the same class of problems considered by Benner and Deni{\ss}en,
we also exploit the fact that changes in $C(v)$ are only low-rank updates,
but we propose a new fast algorithm for efficiently solving such QEPs using a rather different approach.
Our new method also does cubic work once in an offline initialization phase
and only quadratic amount of work in the online phase.
Since many QEPs will typically be solved during the course of optimizing
the viscosities, our approach here can result in a significant speedup
for optimization of viscosities.
At a high level, we propose computing the eigenvalues and eigenvectors of \eqref{eq:qep_eq}
by transforming this QEP into a small sequence of linear eigenvalue problems involving
DPR1 matrices.  By solving these linear subproblems, we can then recover the eigenvalues and
eigenvectors of \eqref{eq:qep_eq}.  Moreover, as these DPR1 matrices are easily converted to DPR1Csym matrices,
we also leverage this special structure in a new fast eigensolver in order to be much more efficient
than standard eigenvalue techniques.

\subsection{Efficient eigenvalue computation for DPR1Csym matrices}
\label{eff_eig_calc}
Let $A \in \C^{2n \times 2n}$ be a DPR1Csym matrix, i.e.,
\begin{equation}
A= D +\rho z z\tp,
\label{A}
\end{equation}
where $D=\diag(d_{1},d_{2},\ldots ,d_{2n}) \in \C^{2n}$ is invertible (so $d_{i}\neq 0$ $\forall i$),
$z=[z_{1} \text{ } z_{2} \text{ }  \cdots \text{ } z_{2n}]\tp \in \C^{2n}$,
and $\rho > 0$.
Note that it is not necessary to consider $\rho \leq 0$,
since if $\rho =0$, then $A$ is diagonal and so obtaining its eigenvalues is trivial,
while if $\rho < 0$, then one can just instead consider $A = - D - \rho zz\tp$.
Furthermore, we assume that
\begin{itemize}
\item $A$ is irreducible, i.e.,
$\forall i,j \in \{1,\ldots, 2n\}$,
$z _{i}\neq 0$ and $d_{i}\neq d_{j}$ if $i \neq j$, and
\item $A$ is diagonalizable.
\end{itemize}
It is unnecessary to consider reducible $A$ matrices, since
$d_i$ is an eigenvalue of $D + \rho z z\tp$, with its corresponding
eigenvector being the $i$th canonical vector,
if and only if $z _{i}=0$ or $d_i = d_j$ for some $j \neq i$ holds
(see, e.g., \cite{XuQ08}).
In other words, such eigenvalues can be easily removed (via exact deflation) to obtain
a smaller DPR1Csym matrix that is irreducible.  Per the following remark,
we will be able to convert eigenvalue problems involving DPR1 matrices into
ones involving DPR1Csym matrices.

 \begin{remark}\label{rem:convert}
 Note that if $D+ \rho u z\tp \in \C^{2n \times 2n}$ is a DPR1 matrix with $u,z \in \C^{2n}$ and $u_i \neq 0$, $z_i \neq 0$ $\forall i$,
 then it can be rewritten as a DPR1Csym matrix with the same eigenvalues.
 Letting
 \begin{equation}
 	\label{matricaS}
	S \coloneqq \diag\left(\sqrt{\tfrac{z_1}{u_1}},\dots,\sqrt{\tfrac{z_{2n}}{u_{2n}}}\right)
	\qquad \text{and} \qquad
	\hat z \coloneqq Su,
\end{equation}
it follows that $u\tp S^2=z\tp$ holds, and so $(\lambda, x)$ is an eigenpair of $D+\rho u z\tp$ if and only if $(\lambda, Sx)$ is an eigenpair of $D+ \rho\hat z\hat z\tp$.
While this transformation requires that $u$ and $z$ only have nonzero entries, this is also easily ensured
via a preprocessing step.
If $u_i = 0$ or $z_i = 0$ for some $i$, then $d_i$ is an eigenvalue of $D$
and it can be removed via exact deflation.
Thus, by first performing a sequence of exact deflations corresponding to the zero entries of $u$ and $z$,
we extract the associated eigenvalues (and eigenvectors) and obtain a smaller DPR1 matrix that can be converted to a DPR1Csym
matrix.
\end{remark}
Thus, with our assumptions above, we need only consider the case of computing
eigenvalues and eigenvectors of DPR1Csym matrices. 

If we were only to consider DPR1 real symmetric matrices, 
then fast standard techniques can be used that exploit the DPR1 structure, e.g., 
divide-and-conquer, where the eigenvalues and eigenvectors
of a tridiagonal matrix are computed by solving a sequence of eigenvalue problems involving DPR1 real symmetric matrices; see \cite{Cup80} and \cite[chapter~5.3.3]{Dem97}.
Of course, the essential properties needed to employ such methods are not present 
for DPR1Csym matrices, the most important being that 
the diagonal elements of $D$ and the eigenvalues of $A$ are no longer interlaced for the complex problem, 
since these values are now in the complex plane as opposed to on the real line.
Thus, we instead consider an approach for DPR1Csym matrices that is inspired by a different approach 
for DPR1 real symmetric matrices~\cite{JakST20}.
The method of~\cite{JakST20} computed eigenpairs using a combination of standard 
and modified Rayleigh quotient iterations (RQI and MRQI, respectively),
but in our setting, the eigenvalues of~\eqref{A} will be complex (and real axis symmetry is not guaranteed),
and we have observed that standard RQI often does not converge.
Moreover, we have also observed that when eigenvalues are close to each other,
the method of \cite{JakST20} often gets stuck oscillating between approximations in such clusters of eigenvalues.
To address these shortcomings, we propose two key modifications, namely, to completely forgo using standard
RQI and to introduce a new dynamic step-size procedure in order to steer our MRQI-based procedure towards
a single eigenvalue in a cluster.
We now present our new method in complete detail.

Since $A$ is complex symmetric and diagonalizable, we have
the following eigendecomposition
\begin{equation}
A=W\Lambda W\tp,
\label{Aeigendec}
\end{equation}
where $\Lambda =\diag(\lambda _{1},\ldots ,\lambda _{2n})$
and $W= \begin{bmatrix} w_{1} & \cdots & w_{2n} \end{bmatrix}$ with $W\tp W = I$
are respectively the eigenvalues and eigenvectors of $A$.
Note that the
eigenvalues of $A$ are the zeros of the secular function (see e.g., \cite{Cup80}):
\begin{equation}
f(\lambda )=1+\rho\sum_{i=1}^{2n}\frac{z _{i}^{2}}{d_{i}-\lambda }=1
+\rho z\tp(D-\lambda I)^{-1}z,  \label{Pick}
\end{equation}
where for $i \in \{1,\ldots,2n\}$,
the eigenvector $w_i$ for eigenvalue $\lambda_i$ is given by
\begin{equation}
w_{i}=\frac{x_{i}}{\left\Vert x_{i}\right\Vert _{2}} \quad\text{with}\quad x_{i}=(
D-\lambda _{i}I) ^{-1}z.  \label{Aeigenvec}
\end{equation}
 The zeros of \eqref{Pick} can be found using different algorithms, e.g., if $A$ is real,
 the eigenvalues can be efficiently and reliably computed via bisection \cite{JakSB15}.
 If $A$ is a DPR1 matrix, one can use, e.g., \verb"mpsolve" from the package MPSolve (see \cite{BinR14}), but this can be costly since \texttt{mpsolve} uses a large amount of extra digits of precision (as opposed to just quad precision).
 If $A$ is a complex symmetric matrix, one can use  MRQI; see \cite{ArbH04,Par74}. 
 Regarding the eigenvector formula given in~\eqref{Aeigenvec},
 this is well known to be numerically unstable, but one option to  work around this problem is to use extended precision; 
 for the DPR1 eigensolver of~\cite{JakSB15}, a fraction of the algorithm
 is implemented in quad precision, and the authors reported that overhead to use this 
 extended precision was very modest, i.e., only 55\% slower than standard double-precision implementations; see \cite[p.~314]{JakSB15}, 

\begin{algorithm}[t]
\caption{Eigensolver for  DPR1 matrices} 
\label{alg:dpr1}
\begin{algorithmic}[1]
\Require
 DPR1 matrix $D + \rho u z\tp$ with $\rho > 0$ and vectors $u, z\in \mathbb{C}^{2n}$ with no zero entries.
\Ensure
	Eigenvalues $\Lambda = \diag(\lambda_1, \ldots, \lambda_{2n})$
	and eigenvectors
		$W= \begin{bmatrix}	w_{1} & \cdots & w_{2n} \end{bmatrix}$.
\State $S \gets $ diagonal matrix from \eqref{matricaS}
\State  $\hat z \gets S u$
\For{$l = 2n,2n-1,2n-2,\ldots,2,1$}
\State $s \gets \mathrm{argmax}_{i \in \{ 1,\ldots m\} }{|d_i|}$
\State $\widehat D \gets D-d_s I$
\State $x \gets e_s$, $\gamma \gets 0$, $\eta \gets 1$ (set initial values)
\While {$|\frac\delta\gamma|<\texttt{tol}$ (stop MRQI once relative change between steps is small)}\label{whilestep}
\If{not converging}
\State $\eta \gets \frac \eta 2$ (reduce $\eta$ and reset other initial values)
\State $x \gets e_s$, $\gamma \gets 0$
\EndIf
\State $\delta \gets \eta x^{T}(\widehat D+\rho \hat z\hat z\tp)x/\|x\|^{2}$ \label{mrqi1}
\State $\gamma \gets \gamma +\delta $
 \State $x \gets (\widehat D-\delta I)^{-1}\hat z$ \label{mrqi2}
 \State $\widehat D \gets \widehat D -\delta I $
\EndWhile
\State $\lambda_l \gets d_s +\gamma $
\State $[D,\hat z]\gets $ Update DPR1Csym matrix by deflating $\lambda_l$ from $D+\rho \hat z \hat z\tp$ via \eqref{deflation_formula}
\EndFor
\State $W\gets \begin{bmatrix} w_{1} & \cdots & w_{2n} \end{bmatrix}$,
		where $w_l$ is computed by~\eqref{Aeigenvec} using $\lambda_l$ 
\State $W \gets S^{-1}W$
\end{algorithmic}
\algcomment{If $u=z$, then $S=I$ in line~1 and so $\hat z = z = u$.
Since $S$ is diagonal, the operations outside of the for loop amount to $\bigO(n^2)$ work, while
each line inside is at most $\bigO(n)$ work.  Thus,
assuming that the number iterations of the while loop is never dependent on $n$,
the total work complexity of \cref{alg:dpr1} is $\bigO(n^2)$.
}	
\end{algorithm}

In our case of $A$ being DPR1Csym, we can consider a modification of the MRQI approach of~\cite{JakST20}
that both additionally takes advantage of its DPR1 structure for efficiency, and introduces our new
step-size procedure to improve the reliability of convergence.
Given a starting $x \in \C^{2n}$, repeat
\begin{equation}\label{mrqi}
\delta \gets \eta\frac{x\tp Ax}{x\tp x}, \quad x \gets (D-\delta I)^{-1}z,
\end{equation}
where $\eta > 0$ is a step size chosen dynamically to enhance convergence to a single eigenvalue. The computation of $x$ comes from the eigenvector formula \eqref{Aeigenvec}.
Once $\delta $ has converged, it can be deflated from $A$ to obtain a new smaller DPR1Csym matrix; see \cite{PanZ11}.
Without loss of generality,  assume eigenvalue $\lambda$ is computed via shift $d_s$ from the diagonal of $D$.
Then deflating $\lambda$ from $A$ results in the DPR1Csym matrix $A_\mathrm{d} \in \C^{2n-1 \times 2n -1}$, where
\begin{subequations}
\label{deflation_formula}
\begin{align}
A_\mathrm d & =D_\mathrm d+\rho z_\mathrm d z\tp_\mathrm d, \quad \text{where} \\
D_\mathrm d &= \diag(d_1,\dots,d_{s-1},d_{s+1},\dots,d_{2n}) \\
(z_\mathrm d)_i &= z_i\,\sqrt{\frac{d_i-d_s}{d_i-\lambda}}, \quad i=1,\dots,s-1,s+1 \dots,2n.
\end{align}
\end{subequations}
The deflation formula comes from shifted inverse power method and Sherman-Morrison-Woodbury (SMW)
formula, and $A$ is always stored implicitly, as two vectors and a scalar.
A detailed pseudocode for our new eigensolver for DPR1 matrices is given in \cref{alg:dpr1}.

\subsection{Efficient eigenvalue computation for QEPs}
\label{dpr1mult}
We now show how to transform our QEP given by \eqref{eq:qep_eq} into multiple connected DPR1 eigenvalue problems.
First, since $M$ and $K$ are symmetric positive definite matrices, there exists a matrix $\Phi$   which simultaneously diagonalizes $M$ and $K$, i.e.,
\begin{equation}
	\label{simult. diag of M,K}
	\Phi\tp K \Phi = \Omega^2= \diag (\omega_1^2, \ldots, \omega_n^2 )
	\quad \text{and} \quad \Phi\tp M \Phi = I,
\end{equation}
where $\omega_1>\cdots>\omega_n>0$ are the undamped frequencies. Moreover, it can be shown that $\Phi$ also diagonalizes
$C_\mathrm{int}$, i.e., $\Phi \tp C_\mathrm{int} \Phi = \alpha \Omega$; for more details,
see  \cite{morBenTT11a,Ves11}.\footnote{While the eigensolver we propose in this section assumes that $\Phi$ diagonalizes $C_\mathrm{int}$, note that our choice to use critical damping, i.e., \eqref{damp_int}, is not required. In particular, our approach can be applied to any internal damping that corresponds to a modally damped system, which is a usual assumption when vibrational mechanical systems are considered.}
Thus, we can linearize the QEP given in \eqref{eq:qep_eq} to obtain the standard eigenvalue problem
\begin{subequations}
\begin{align}
	\label{eq:EP}
	A(v) y(v) &= \lambda(v)  y(v), \quad \text{where} \\
	\label{A form}
	A(v) &= 	
	\begin{bmatrix}
		0 & \Omega \\
	  	-\Omega & -\Phi\tp C(v)\Phi
	 \end{bmatrix}
	 =
	 \begin{bmatrix}
	 	0 & \Omega \\
		  -\Omega & -\alpha\Omega
	\end{bmatrix}
	-
	\begin{bmatrix}
	  	0  \\
		\Phi\tp G
	\end{bmatrix}
	\begin{bsmallmatrix}
		v_1 &  &  \\
	   	& \ddots &  \\
	  	&  & v_r
	\end{bsmallmatrix}
	\begin{bmatrix}  0 &  G\tp \Phi \end{bmatrix}, \\
	y(v) &= \begin{bmatrix}\Omega \Phi^{-1} x(v)\\
\lambda(v)\Phi^{-1}x(v)
\end{bmatrix}.
\end{align}
\end{subequations}
Let $P \in \R^{2n \times 2n}$ be the perfect shuffle permutation, which splits a set of even cardinality into two sets of equal cardinality and interleaves them,
i.e., it maps the $k$th entry as follows:
\[
	k \mapsto
	\begin{cases}
      		2k-1, & \mbox{if } k\leq n \\
                 2(k-n), & \mbox{if } k>n.
 	\end{cases}
\]
Now dropping the dependency on $v$ for brevity and using $PP\tp = I$ and $\widehat A = P\tp A P$,
multiplying \eqref{eq:EP} on the left by $P\tp$ yields the eigenvalue problem
\begin{subequations}
\begin{align}
	\label{eq:blockdiag}
	\widehat A P\tp y
  	&= \left(
	\begin{bsmallmatrix}
          	D_1 & & \\
              	& \ddots &\\
           	& &D_n\\
        	\end{bsmallmatrix}
	-
	\widehat G
	\begin{bsmallmatrix}
		v_1 &  &  \\
	   	& \ddots &  \\
	  	&  & v_r
	\end{bsmallmatrix}
	\widehat G\tp
	\right)
	P\tp y = \lambda P\tp y,  \quad \text{where} \\
	\label{eq:Di_gj}
	D_i &= \begin{bmatrix}
            	0 & \omega_i \\
           	-\omega_i & -\alpha \omega_i
        	\end{bmatrix}
	\quad \text{and} \quad
	\widehat G =
	P\tp
	\begin{bmatrix}
	  	0  \\
		\Phi\tp G
	\end{bmatrix}.
\end{align}
\end{subequations}
Let $\Psi_i$ be the matrix which diagonalizes matrix $D_i$ and consider the matrices
\begin{equation}
	\label{psi_mat}
  	\Psi = \begin{bsmallmatrix}
        		\Psi_1 &  &  \\
              	& \ddots &  \\
              	&  & \Psi_n
     	\end{bsmallmatrix},
	\quad
	D = \Psi^{-1}
		\begin{bsmallmatrix}
          		D_1 & & \\
              		& \ddots &\\
 		        & &D_n\\
        		\end{bsmallmatrix}
		\Psi,
	\quad
	U = \Psi^{-1}\widehat G,
	\quad \text{and} \quad
	Z = \Psi\tp \widehat G,
\end{equation}
noticing that $\Psi$ diagonalizes the block diagonal matrix from \eqref{eq:blockdiag} into $D$.
Thus, considering $\widetilde A = \Psi^{-1} \widehat A \Psi$
and multiplying \eqref{eq:blockdiag} on the left by $\Psi^{-1}$,
we further transform the eigenvalue problem into one involving a diagonal matrix plus a low-rank update
\begin{subequations}
\begin{align}
	\label{tildeA}
	\widetilde Aw &= \lambda w,
	\quad \text{where} \\
	\widetilde A &= D - U
		\begin{bsmallmatrix}
		v_1 &  &  \\
	   	& \ddots &  \\
	  	&  & v_r
		\end{bsmallmatrix}Z\tp
	= D -\sum\limits_{j=1}^r v_j u_jz_j\tp,
	\quad
	w = \Psi^{-1}P\tp y,
\end{align}
\end{subequations}
and $u_j$ and $z_j$ are respectively the $j$th columns of $U$ and $Z$.
Since matrices $\Phi$, $\Omega$, $P$, $D$, $U$, and $Z$ are all independent of $v$,
obtaining the low-rank structure of $\widetilde A$ can be precomputed once in an offline  process.
In fact, $\Phi$, $\Omega$, $P$, and $D$ are also independent of the damping positions
specified by the matrix $G$, and so, when optimizing viscosities for multiple configurations of damping positions,
these matrices need only be calculated once, while
computing $U$ and~$Z$ for each configuration is cheap.

\begin{algorithm}[t]
\caption{QEP eigensolver for \eqref{eq:qep_eq}}
\label{alg:dpr1mult}
\begin{algorithmic}[1]
\Require $M$ and $K$ from \eqref{eq:qep_eq}, $\alpha \geq 0$ for $C_\mathrm{int}$ from \eqref{damp_int},
	$\Phi$ and $\Omega$ from \eqref{simult. diag of M,K},
	$\Psi$, $D$, $U$, and $Z$ from \eqref{psi_mat},
	and $v \in \R^r$.
 \Ensure Eigenvalues $L$ and eigenvectors $X$ of QEP \eqref{eq:qep_eq}
 	\State $L_0\gets D$
 	\For{\texttt{$j=1,\dots, r$}}
      	 	\State $[\tilde u_j,\tilde z_j] \gets$ $j$th columns of $U$ and $Z$, respectively
        		\State $[L_j,\xi_j] \gets$ eigenvalues and eigenvectors of
				$L_{j-1} - v_j \tilde u_j \tilde z_j\tp$ computed
			 	by \cref{alg:dpr1}
		\State $ U \gets \xi_j^{-1}U$
        		\State $ Z \gets \xi_j\tp Z$
      	\EndFor
      	\State $L \gets L_r$
      	\State $\Xi \gets P \Psi \xi_1 \cdots \xi_r$, where $P \in \R^{2n \times 2n}$ is the perfect shuffle permutation
      	\State $X \gets \Phi \Omega^{-1}\Xi(\texttt{$1$:$n$,:})$ (Take the first $n$ rows of $\Xi$)
	\State $X \gets$ each column (an eigenvector) of $X$ gets refined according to \cref{rem:inv_smw}
\end{algorithmic}
\algcomment{For simplicity of the pseudocodes in this paper,
we assume that vectors $\tilde u_j$ and $\tilde z_j$ never have zero entries, scalar $v_j \neq 0$,
and $L_{j-1} - v_j \tilde u_j \tilde z_j\tp$ is actually given to \cref{alg:dpr1} as
$L_{j-1} + |v_j| (-\sgn(v_j)\tilde u_j) \tilde z_j\tp$ to adhere to its convention that $\rho > 0$.
If $v_j=0$, obtaining the eigenvalues and eigenvectors is immediate (so \cref{alg:dpr1} is not needed),
while if either $\tilde u_j$ or $\tilde z_j$ contain zero entries, then per \cref{rem:convert},
exact deflation is first used to remove the
corresponding eigenvalues, which are on the diagonal of~$L_{j-1}$,
and then \cref{alg:dpr1} is called on the resulting deflated DPR1 matrix
to obtain the remaining eigenvalues and eigenvectors.
Finally, note that by design of \cref{alg:dpr1,alg:dpr1mult},
for $j > 1$ in line~4, \cref{alg:dpr1} is warm started by using the eigenvalues
of the previous DPR1 eigenvalue problem as shifts for the next DPR1 eigenvalue problem.
}		
\end{algorithm}

We now show how \cref{alg:dpr1} can be iteratively applied to portions of $\widetilde A$
in order to recover all the eigenvalues and eigenvectors of \eqref{tildeA}.
Let $\widetilde A_1 = D - v_1 u_1 z_1\tp$ and suppose it is diagonalizable, i.e.,
it has eigendecomposition $\widetilde A_1 = \xi_1 L_1 \xi_1^{-1}$, where matrices $\xi_1$ and $L_1$ respectively
contain the eigenvectors and eigenvalues of $\widetilde A_1$.
Then multiplying \eqref{tildeA} on the left by $\xi_1^{-1}$ and separating out $\widetilde A_1$, we obtain
the transformed eigenvalue problem
\begin{equation}
	\xi_1^{-1}\left(\widetilde A_1 - \sum_{j=2}^r v_j u_j z_j\tp \right) w
	=\left(L_1 - \sum_{j=2}^r v_j \xi_1^{-1} u_j z_j\tp \xi_1 \right) \xi_1^{-1} w
	= \lambda \xi_1^{-1} w.
\end{equation}
If $\widetilde A_2 = L_1 - v_2 \tilde u_2 \tilde z_2\tp$ is also diagonalizable,
where $\tilde u_2 = \xi_1^{-1}u_2$ and $\tilde z_2 = \xi_1\tp z_2$,
we can again similarly transform the eigenvalue problem via the eigendecomposition
$\widetilde A_2 = \xi_2 L_2 \xi_2^{-1}$.  We keep applying these transformations for $j=1,\ldots,r$ by computing
the eigendecompositions
\begin{subequations}
\begin{align}
	\widetilde A_j = \xi_j L_j \xi_j^{-1},
	\qquad \text{where} \qquad
	\widetilde A_j &= L_{j-1} - v_j \tilde u_j \tilde z_j\tp, & L_0 &= D, \\
	\tilde u_j &= \xi_{j-1}^{-1} \cdots \xi_1^{-1} u_j, & \tilde u_1 &= u_1, \\
	\tilde z_j &= \xi_{j-1}\tp \cdots \xi_1\tp z_j, & \tilde z_1 &= z_1.
\end{align}
\end{subequations}
Assuming all the $\widetilde A_j$ matrices are indeed diagonalizable,
we finally obtain
\begin{equation}
 	L_r \left( \xi_r^{-1} \cdots \xi_1^{-1} w \right) = \lambda \left(\xi_r^{-1} \cdots \xi_1^{-1} w\right),
\end{equation}
and so we have recovered the eigenvalues of \eqref{tildeA} and can reconstruct its eigenvectors as well.

\begin{remark}
\label{rem:inv_smw}
As a final step of our algorithm, note that we also follow a suggestion of Taslaman~\cite[section~4.3]{Tas15}
to refine the accuracy of computed eigenvectors by doing a single step of inverse iteration for each eigenvector;
due to the particular structure of $C(v)$, the SMW formula can be used to do 
a single-step of inverse iteration in only $\bigO(n)$ work per eigenvector.
Similar application of the SMW formula in damped systems for efficient computations can be found in~\cite{TruV09,TruTV15,KuzTT16,morTomBG18,morBeaGT20}.
\end{remark}

Pseudocode for our complete QEP eigensolver is given in~\cref{alg:dpr1mult}.
We note that assuming that all matrices $\widetilde A_j$ are diagonalizable is quite standard
(see, e.g., \cite{morBenKTetal16,JakST20}), and we have not observed any issues in practice
with this assumption.  

We now turn to the work complexity of~\cref{alg:dpr1mult}.
Recall that the work complexity of~\cref{alg:dpr1} is~$\bigO(n^2)$,
and since we assume that the number of dampers $r$ is small, i.e., $r \ll n$,
we will treat $r$ as a constant.
Inside the for loop of~\cref{alg:dpr1mult}, lines~5 and 6 are potentially more than $\bigO(n^2)$ work using
standard techniques, but since the $\xi_j$'s are Cauchy-like matrices,
matrix-vector multiplication and linear solves can be done in approximately linear time,
and so the total cost of the loop remains at $\bigO(n^2)$.
Meanwhile, forming~$\Xi$ in line~9 is also not more than $\bigO(n^2)$ work,
since $P$ and $\Psi$ are sparse matrices and $\xi_1 \cdots \xi_r$ is a product of linked Cauchy-like matrices,
which can also be done in (approximately) linear time;
for full details, see \cite{PanZ00,JakST20}.
In line~10, applying the inverse of $\Omega$ is cheap (recall that $\Omega$ is a diagonal matrix),
while if $\Phi$ is a sparse matrix, then
obtaining the full matrix of eigenvectors~$X$ is also at most $\bigO(n^2)$ work.
However, note that while evaluating the functions in~\eqref{eq:model1} and \eqref{eq:model2}
requires that we obtain all of the eigenvalues of \eqref{eq:qep_eq},
only a handful of the corresponding eigenvectors are needed to compute the corresponding gradients.
For our setting of optimizing viscosities using gradients, in line~10, we can
selectively compute the handful of relevant of eigenvectors, i.e.,
we apply~$\Phi$ to the few corresponding columns of~$\Omega^{-1}\Xi(\texttt{$1$:$n$,:})$
in order to not exceed $\bigO(n^2)$ work when~$\Phi$ is dense.
Finally, per \cref{rem:inv_smw}, refining the eigenvectors using inverse iteration is only $\bigO(n)$ work
per eigenvector.
Hence, the overall work complexity of \cref{alg:dpr1mult} is quadratic.

In terms of constant factors, the total cost of \cref{alg:dpr1mult} is dominated by line~4, i.e.,
the~$s$~calls to \cref{alg:dpr1}.
As such, provided $n$ is small enough,
it is not always critical to implement lines~5, 6, 9, and 10 in \cref{alg:dpr1mult} as discussed above in order to attain
the theoretical work complexity result (but those steps should indeed be taken if $n$ is very large).
Also, we note that the subsequent calls to \cref{alg:dpr1} can be faster than the first one because, by our design of the two algorithms,
we are actually warm starting \cref{alg:dpr1} by choosing the initial shifts as the eigenvalues of the previous eigenvalue problem.
Thus, when two consecutive DPR1 eigenvalue problems have quite similar spectra, which is not uncommon,
we end up having excellent initial shift choices for which
to accelerate the convergence of~\cref{alg:dpr1} in line~4 of~\cref{alg:dpr1mult}.

\section{The frequency-weighted damping optimization algorithm}
\label{sec:damp_opt}
We are now ready to present our new algorithm for frequency-weighted damping of QEPs.
We begin with the offline phase, which simply precomputes the matrices from \eqref{simult. diag of M,K} and~\eqref{psi_mat}
so that \cref{alg:dpr1mult} can be used to evaluate all the eigenvalue-based functions (and their gradients)
that appear in \eqref{eq:model1} and \eqref{eq:model2}.
While this offline phase has a $\bigO(n^3)$ work complexity, it only needs to be done once.
We emphasize again that the simultaneous diagonalization part of the offline phase is also independent of the damping positions,
and so it only needs to be performed once for all different damping positions.

For the online phase, by using \cref{alg:dpr1mult}, evaluating all the functions (and their gradients
described in \cref{sec:opt_solve})
in \eqref{eq:model1} or \eqref{eq:model2}
for a given vector $v \in \R^r$ is then only $\bigO(n^2)$~work,
as opposed to $\bigO(n^3)$ via standard eigensolvers.
In terms of the overall cost, this is a significant savings
as we expect to require many function evaluations before converging to a stationary point
of \eqref{eq:model1} or \eqref{eq:model2}, particularly since these are nonsmooth optimization problems.
To find solutions of \eqref{eq:model1} and \eqref{eq:model2}, we use GRANSO;
a high-level description of our method is given in \cref{alg:opt}.

\begin{algorithm}[t]
\caption{Frequency-weighted damping optimization algorithm}
\label{alg:opt}
\begin{algorithmic}[1]
\Require $M$ and $K$ from \eqref{eq:qep_eq}, $\alpha \geq 0$ for $C_\mathrm{int}$ and $G$ from \eqref{damp_matr},
	set of $k$ ellipses $\mathcal{E}$, weights~$[\phi_1, \ldots, \phi_k]$ with each $\phi_j \in (0,1]$ for ellipse $E_j \in \mathcal{E}$,
	$\eta \geq 0$, $\texttt{tol}_\mathrm{sa} < 0$, initial viscosity values $v_\mathrm{init} \in \nnR^r$,
	and $\texttt{approach} \in \{1,2\}$.
\Ensure Computed for optimized viscosities $v_\mathrm{opt} \in \nnR^r$ for either \eqref{eq:model1} or \eqref{eq:model2}
\State {\bf Offline stage:} (Set up for computing functions and gradients via \cref{alg:dpr1mult})
  \State $[\Phi,\Omega] \gets$ matrices from \eqref{simult. diag of M,K} (Diagonalization)
  \State $[\Psi,D,U,Z] \gets$ matrices from \eqref{psi_mat} (Linearize and construct low-rank structure)
  \State {\bf Online stage:} (Optimize viscosities using GRANSO and \cref{alg:dpr1mult})
 \If { $\texttt{approach} = 1$ }
 	\State $v_\mathrm{opt} \gets$ solution returned by GRANSO for \eqref{eq:model1} initialized at $v_\mathrm{init}$
\Else
	\State $v_\mathrm{opt} \gets$ solution returned by GRANSO for \eqref{eq:model2} initialized at $v_\mathrm{init}$
\EndIf
\end{algorithmic}
\end{algorithm}

\section{Numerical experiments}\label{nume_ex}
\begin{figure}[t]
\begin{center}
\centering
%
%
%

 \begin{tikzpicture}[scale=1, every node/.style={scale=0.8}]
\tikzstyle{spring}=[thick,decorate,decoration={zigzag,pre
length=0.3cm,post length=0.3cm,segment length=6}]
\tikzstyle{damper}=[thick,decoration={markings,
  mark connection node=dmp,
  mark=at position 0.5 with
  {
    \node (dmp) [thick,inner sep=0pt,transform
    shape,rotate=-90,minimum width=10pt,minimum height=2pt,draw=none]
    {};
    \draw [thick] ($(dmp.north east)+(1.5pt,0)$) -- (dmp.south east)
    -- (dmp.south west) -- ($(dmp.north west)+(1.5pt,0)$);
    \draw [thick] ($(dmp.north)+(0,-3pt)$) -- ($(dmp.north)+(0,3pt)$);
  }
}, decorate]
\tikzstyle{springdot}=[thick,decoration={markings,
  mark connection node=sdt,
  mark=at position 0.5 with
  {
  \node (sdt) [thick,inner sep=0pt,transform shape,rotate=-90,minimum
  width=0.85cm,minimum height=0.50cm,draw=none] {};
      \draw [dotted, thick, color=blue] ($(sdt.north)$) --
      ($(sdt.south)$);
  }
}, decorate]

\tikzstyle{ground}=[fill,pattern=north east lines,draw=none,minimum
width=0.75cm,minimum height=0.3cm]

\newcommand{\drawLinewithBG}[2]
{
    \draw[white,myBG]  (#1) -- (#2);
    \draw[black,very thick] (#1) -- (#2);
}

\tikzstyle myBG=[line width=3pt,opacity=1.0]

\node (M) [draw,outer sep=0pt,thick,minimum width=1.2cm, minimum
height=1.2cm,color=red, fill=red!20!white] at (0,0) {$\color{black}
m_1$};
\node (M2) [draw,outer sep=0pt,thick,minimum width=1.2cm, minimum
height=1.2cm,color=red,fill=red!20!white] at (2,0) {$\color{black}
m_2$};
\node (M3) [draw,outer sep=0pt,thick,minimum width=1.2cm, minimum
height=1.2cm,color=red,fill=red!20!white] at (4,0) {$\color{black}
m_{n-1}$};
\node (M4) [draw,outer sep=0pt,thick,minimum width=1.2cm, minimum
height=1.2cm,color=red,fill=red!20!white] at (6,0) {$\color{black}
m_n$};

\node (LW)
[ground,rotate=-90,anchor=north,xshift=-1cm,yshift=-2.0cm,minimum
width=3cm, style={draw,outer sep=0pt,thick}] at (M.south) {};

\node (RW)
[ground,rotate=90,anchor=north,xshift=1cm,yshift=-2.5cm,minimum
width=3cm, style={draw,outer sep=0pt,thick}] at (M4.south west) {};

\draw [spring,color=blue] ({-1.60,0}) -- ({-0.48,0});
\draw [spring,color=blue] (M2.180) -- ($(M.north
east)!(M.180)!(M.south east)$);
\draw [damper] (M2.150) -- ($(M.north east)!(M.150)!(M.south east)$);

\draw [springdot,color=blue] (M3.180) -- ($(M2.north
east)!(M2.180)!(M2.south east)$);
\draw [spring,color=blue] (M4.180) -- ($(M3.north
east)!(M3.180)!(M3.south east)$);
\draw [spring,color=blue] ({6.48,0}) -- ({7.52,0});
\draw [damper] ({0,1.1}) -- ({-1.1,1.1});
\draw [thick] ({0,1.1}) -- ({0,.5});

\draw [damper] ({4,1.1}) -- ({2.9,1.1});
\draw [thick] ({4,1.1}) -- ({4,.5});

\node at (-.8,-.3) {$k_1$};
\node at (-.7,.7) {$v_1$};

\node at (1,-.3) {$k_2$};
\node at (1,.7) {$v_2$};

\node at (3.3,.7) {$v_3$};
\node at (5,-.3) {$k_n$};
\node at (7.15,-.3) {$k_{n+1}$};

\node (LW)
[ground,rotate=-90,anchor=north,xshift=-1cm,yshift=-2.0cm,minimum
width=1cm, style={draw,outer sep=0pt,thick}] at (.5,.3) {};
\node (LW)
[ground,rotate=-90,anchor=north,xshift=-1cm,yshift=-2.0cm,minimum
width=1cm, style={draw,outer sep=0pt,thick}] at (4.5,.3) {};

\end{tikzpicture}

%
%
\caption{Diagram of an $n$-mass oscillator.}
\label{mck_fig}
\end{center}
\end{figure}
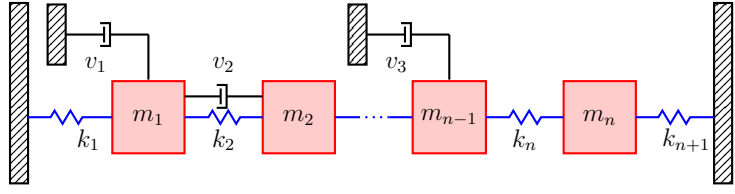

All experiments were done in \MATLAB\ R2021a using a mid-2020 13" MacBook Pro with an Intel Core i5-1038NG7 CPU (quad core) and 16GB of RAM
 running macOS 10.15.7.
 Our code for replicating all experiments reported here is provided in the supplementary material.
For the values of $n$ in our experiments here,
it sufficed to implement lines~5, 9, and 10 of~\cref{alg:dpr1mult} using standard techniques
and compute all the eigenvectors, as opposed to leveraging the Cauchy-like structure
and possibly selectively computing eigenvectors.
 As test problems, we used various instances of an $n$-mass oscillator; see \cref{mck_fig}.
For this mechanical system, we have the following matrices
\begin{subequations}
\label{eq:nmass}
\setlength\arraycolsep{4pt}
\begin{align}
	M &= \diag(m_1,m_2,\dots,m_n), \\
    	K &=
	\begin{bsmallmatrix}
	k_1+k_2	& -k_2  	& 			& 			\\
        	-k_2   	& \ddots 	& \ddots   		& 			\\
          		&  \ddots 	& \ddots  		& -k_{n}  		\\
               		&      		&  -k_n 		& k_n+k_{n+1}
        	\end{bsmallmatrix}, \\
	C_\mathrm{ext}(v) &= v_1e_je_j\tp+v_2(e_k-e_{k+1})(e_k-e_{k+1})\tp+v_3e_le_l\tp,
\end{align}
\end{subequations}
where $e_j$ denotes the $j$th canonical vector, and $v_1,v_2,v_3 \geq 0$ are the viscosity values.
 In $C_\mathrm{ext}(v)$, the $e_je_j\tp$ and $e_le_l\tp$ terms respectively mean
that there are grounded dampers on masses $m_k$ and $m_l$, while $(e_k-e_{k+1})(e_k-e_{k+1})\tp$
indicates that masses $m_k$ and $m_{k+1}$ are connected by a damper. Thus, for \cref{mck_fig},
$C_\mathrm{ext}(v)$ is defined using $j=k=1$ and $l=n-1$.
Considering \eqref{damp_matr},
we also have that $C_\mathrm{ext}(v) = G\, \diag(v_1,v_2,v_3) G\tp$,
where $G = \begin{bsmallmatrix} e_j & e_k-e_{k+1} & e_l\end{bsmallmatrix}$.

\subsection{Validating \cref{alg:dpr1mult}}
\label{examp1}
To assess the efficiency and accuracy of our new eigensolver (\cref{alg:dpr1mult})
for solving \eqref{eq:qep_eq},
we used instances of \eqref{eq:nmass} with orders $n=200,400,600,\ldots,2000$.
For each value of $n$, we defined matrix $M$ using $m_i=10+990(\frac{i - 1}{n-1})$ for $i=1,\ldots,n$
and matrix $K$ using $k_i=5$ for $i=1,\ldots,n+1$,
and created two problems with different configurations of dampers,
by defining two \mbox{$C(v) = C_\mathrm{int} + C_\mathrm{ext}(v)$} matrices.
We used $\alpha=0.004$ to define $C_\mathrm{int}$, while
the two versions of $C_\mathrm{ext}(v)$ were defined via choosing $j$, $k$, and $l$
as follows:
\[
	\label{conf_set}
	\text{Config A: } (j,k,l)=\left(\tfrac{n}{10},\tfrac{3n}{10},\tfrac{5n}{10}\right)
	\quad \text{and} \quad
	\text{Config B: } (j,k,l)=\left(\tfrac{3n}{10},\tfrac{7n}{10},\tfrac{9n}{10}\right).
\]
Using randomly generated viscosity values for each $n$, specifically $v = 0.1 + \texttt{rand($3$,$n$)}$,
we solved the resulting QEPs with \cref{alg:dpr1mult} and other solvers for comparison purposes.
For direct QEP solvers, we tested \texttt{polyeig} and \texttt{quadeig}.
We also benchmarked \cref{alg:dpr1mult} against a much simpler version of our algorithm,
which also first computes matrix~$\Phi$ to simultaneously diagonalize $M$ and $K$, per \eqref{simult. diag of M,K},
but then forgoes taking any advantage of low-rank structure and instead just computes the eigenvalues of $A(v)$ defined
in \eqref{A form} via calling \texttt{eig} on this standard eigenvalue problem;
we refer to this simpler method as \texttt{eig($A$)}.

\begin{figure}[t]
\centering
\subfloat[Config A]{
\resizebox*{\figsizes}{!}{\includegraphics[trim=0cm 0cm 0cm 0cm,clip]{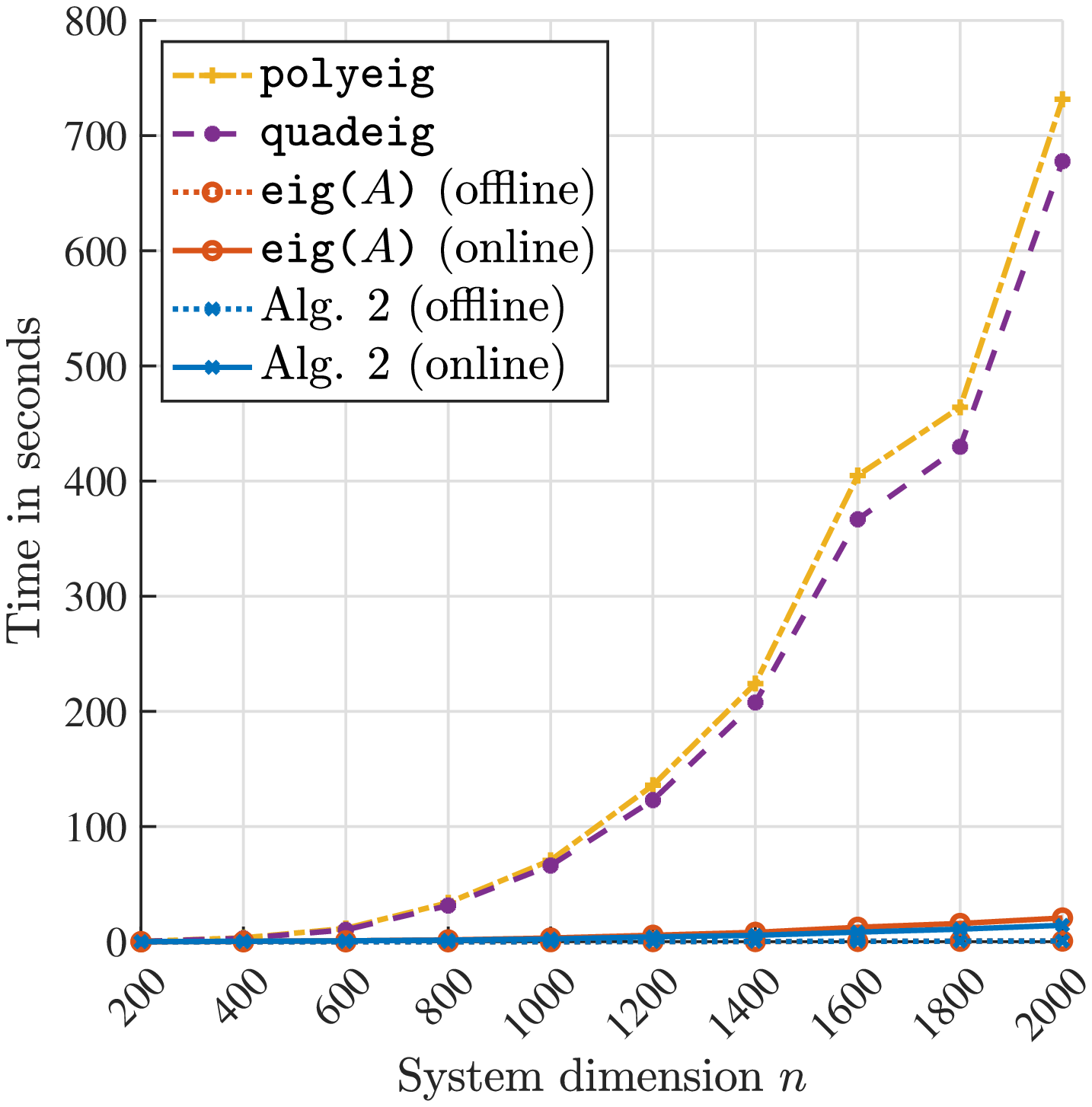}}
}
\subfloat[Zoomed in view]{
\resizebox*{\figsizes}{!}{\includegraphics[trim=0cm 0cm 0cm 0cm,clip]{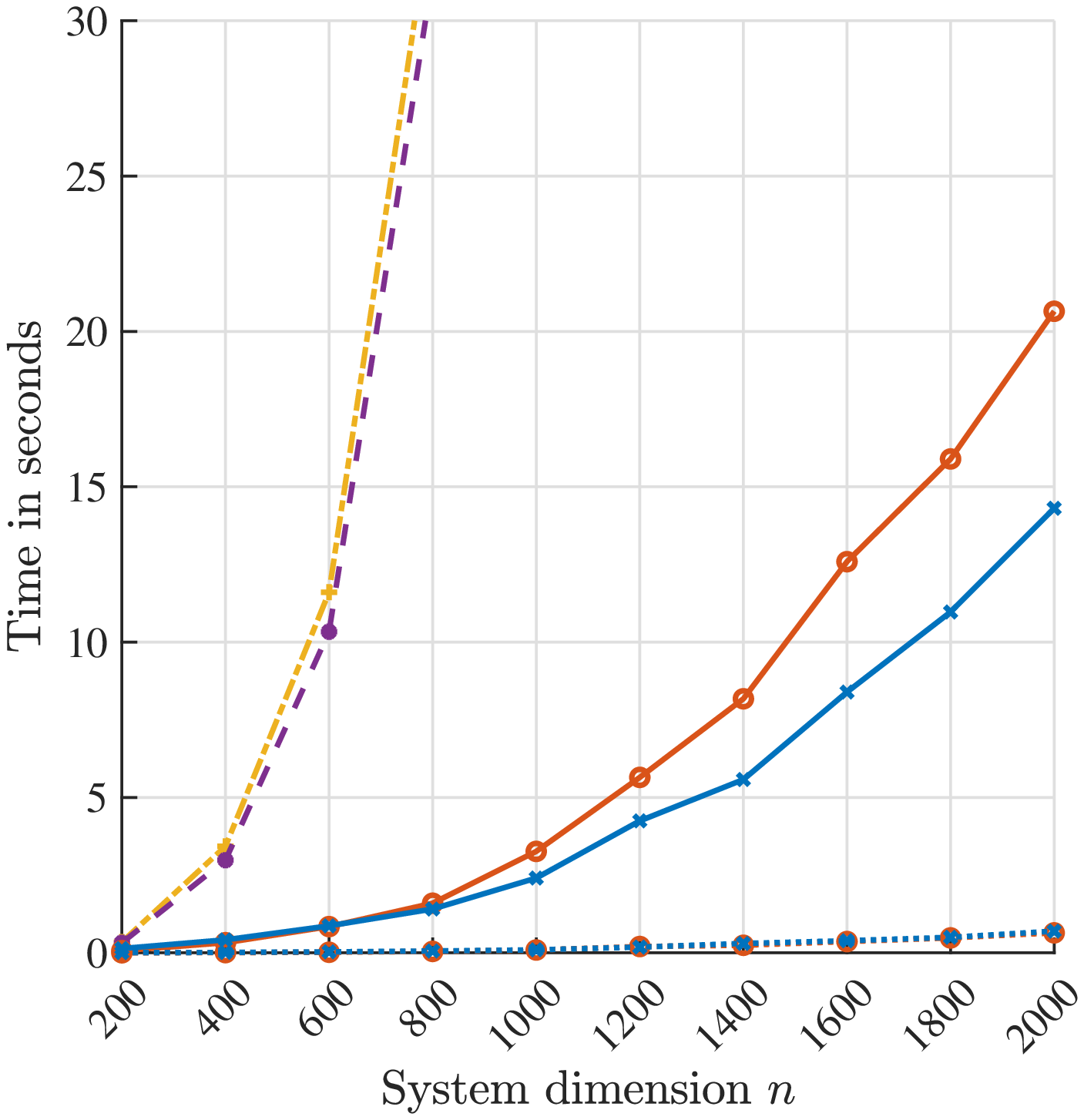}}
\label{fig:scaling_zoom}
}
\caption{The overall running times as the system dimension $n$ is increased.
For the direct solvers, \texttt{polyeig} and \texttt{quadeig}, the overall running times
to solve instances of \eqref{eq:qep_eq} are shown.
For the indirect solvers, \texttt{eig($A$)} and \cref{alg:dpr1mult},
the overall running times are separated into their offline and online parts,
where the offline cost for both is computing $\Phi$ from \eqref{simult. diag of M,K}.
}
\label{fig:scaling}
\end{figure}

In~\cref{fig:scaling}, we show the overall running times, recorded using \texttt{tic} and \texttt{toc},
for all the different eigensolvers as $n$ increases for Config~A.  As can be seen,
the costs of the direct solvers, \texttt{polyeig} and \texttt{quadeig}, quickly increase as $n$ does.
Meanwhile, \texttt{eig($A$)} and \cref{alg:dpr1mult} are much faster, with their respective costs
also increasing at a much slower rate with respect to~$n$.  Indeed,
already by $n=400$, \texttt{eig($A$)} and \cref{alg:dpr1mult} are about an order of magnitude faster than the direct solvers.
Moreover, for $n=2000$, \texttt{eig($A$)} is about 35~times faster than \texttt{polyeig},
while~\cref{alg:dpr1mult} is 51~times faster than \texttt{polyeig}.
Comparing \texttt{eig($A$)} and \cref{alg:dpr1mult} to each other (see \cref{fig:scaling_zoom}),
we see the cost of the latter indeed grows more slowly with respect to $n$, 
and that hidden constant term in the work complexity for \cref{alg:dpr1mult}
is not an issue for overall efficiency in practice.
We note that an implementation of~\cref{alg:dpr1mult} in a compiled language
and that takes advantage of the Cauchy-like structure of the $\xi_j$ matrices should
be many times faster than our prototype implementation that we have used here,
which recall, is coded in \MATLAB\ and does not yet take advantage of Cauchy-like structure.
We also performed the same scaling experiment for Config~B, which resulted
in plots very similar to those shown in~\cref{fig:scaling}; as such, we omit these additional plots here.

\begin{figure}[t]
\centering
\subfloat[Config A]{
\resizebox*{\figsizes}{!}{\includegraphics[trim=0cm 0cm 0cm 0cm,clip]{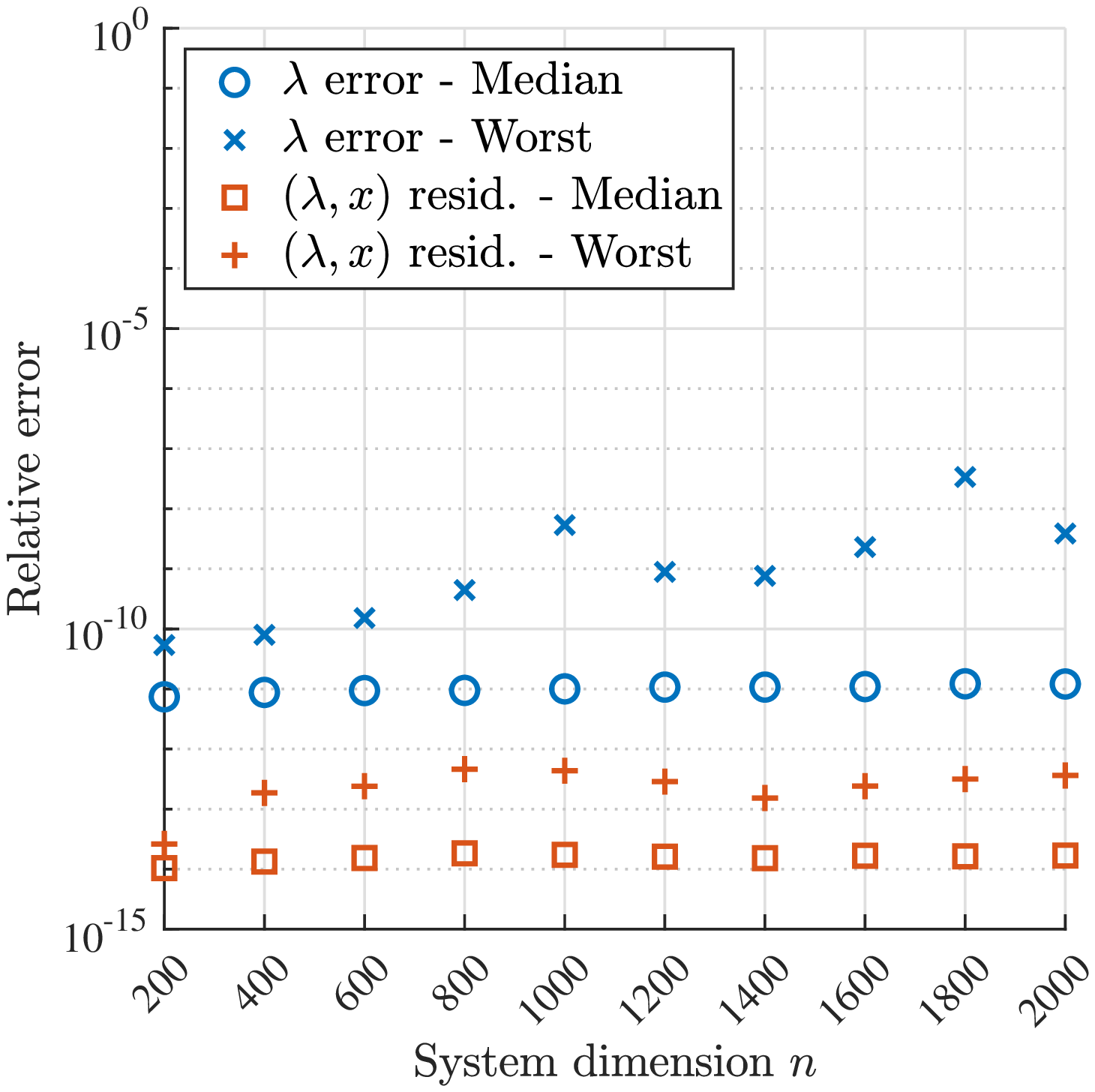}}
}
\subfloat[Config B]{
\resizebox*{\figsizes}{!}{\includegraphics[trim=0cm 0cm 0cm 0cm,clip]{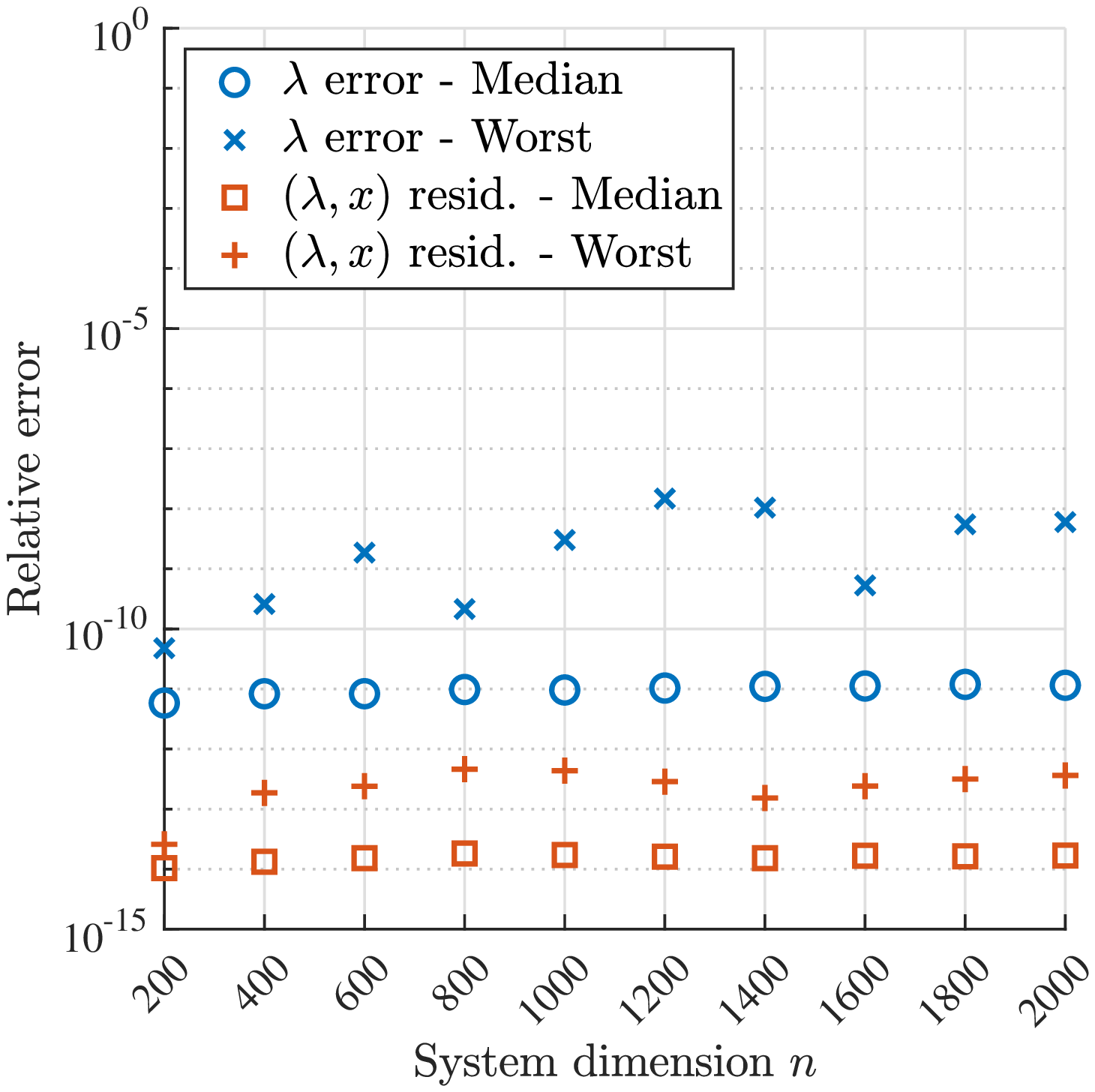}}
}
\caption{The median and worst relative errors of eigenvalues computed by~\cref{alg:dpr1mult} with
respect to the eigenvalues computed by~\texttt{polyeig} are denoted by the ``$\lambda$ error" markers,
while the median and worst eigenpair residuals are denoted by ``$(\lambda,x)$ resid." markers.}
\label{fig:rel_err}
\end{figure}

In order to show the accuracy of \cref{alg:dpr1mult},
we compared its computed eigenvalues with those computed by \texttt{polyeig},
and for each computed eigenvalue-eigenvector pair $(\lambda,x)$, we computed the spectral 
norm of~\eqref{eq:qep_eq} with this eigenpair plugged into it.
For each problem, we paired the two sets of computed eigenvalues greedily, i.e.,
by taking the closest pair of values across the two sets, removing this ``match", and
then repeating this procedure until all computed eigenvalues were paired.
For each matched pair of eigenvalues, we computed the relative errors in the real and imaginary parts separately,
which we denote $\delta_\Re$ and $\delta_\Im$, and then used $\max \{ |\delta_\Re|, |\delta_\Im|\}$
as an overall measure of the error in the computed pair.
Then for each problem, we computed the median and worst (largest) of these error measurements over
the entire computed spectrum.  
Similarly, we computed the median and worst (largest) errors of the norm of~\eqref{eq:qep_eq}
over all the computed eigenpairs~$(\lambda,x)$.
We show the resulting error measurements in \cref{fig:rel_err}
for both Config~A and Config~B across all values of~$n$ tested.
As can be seen, the results are essentially the same for both configurations.
For comparing the accuracy of the eigenvalues with respect to those computed by \texttt{polyeig},  
the median error was always about $10^{-11}$, while the worst error rose from about
$10^{-10}$ to a bit over~$10^{-8}$ as $n$ increased from~200 to~2000;
we saw very similar eigenvalue errors when comparing \cref{alg:dpr1mult} to \texttt{quadeig}
and even when comparing \texttt{polyeig} to \texttt{quadeig}.
Meanwhile, the eigenpairs residuals were in the worst case still under $10^{-12}$ 
with the median error being about~$10^{-14}$, thus demonstrating
that \cref{alg:dpr1mult} is indeed computing eigenvalues and eigenvectors to good accuracy.
Note that our \MATLAB\ implementation of~\cref{alg:dpr1mult} only uses double precision
 and that implementing the key parts of \cref{alg:dpr1} using quad precision should improve 
 the accuracy of~\cref{alg:dpr1mult}; in this case, \cref{alg:dpr1} and \cref{alg:dpr1mult} would be 
mixed-precision implementations.

\subsection{Validating \cref{alg:opt} for Approaches~1 and 2}
\label{examp3}
To demonstrate our new approaches for optimizing viscosities via nonsmooth constrained optimization
and our new models, Approach~1 (Fixed ellipses) and Approach~2 (Variable ellipses),
we used additional instances of the three-damper $n$-mass oscillator defined by the matrices in \eqref{eq:nmass}.
For these experiments, we used $n=1000$ and defined $M$ and $K$ via respectively
setting \mbox{$m_i=m_{n+1-i}=\frac{2n-i}{200}$} for $i=1,\ldots,\frac n 2$
and $k_i=5$ for $i=1,\ldots,n+1$.
To define $C(v)$, we used $(j,k,l) = (100,400,900)$ to specify the configuration of dampers
given by matrix $G$ in $C_\mathrm{ext}(v)$
and used various values of $\alpha$ (to be reported momentarily) for $C_\mathrm{int}$.

For the online optimization phase of~\cref{alg:opt},
we used GRANSO's default parameters except we set \texttt{opts.maxit=100},
always initialized GRANSO from \mbox{$v_\mathrm{init} = \texttt{ones(3,1)}$},
and set \texttt{opts.mu0=10000}.  This last change, which multiplies the objective function by 10000,
was simply done for rescaling reasons, i.e.,
so that the value of the objective function at $v_\mathrm{init}$ was about one for all of our test problems;
 in practice, \texttt{opts.mu0} can be easily determined from the specific problem or one can use GRANSO's
 automatic pre-scaling feature.
Since \eqref{eq:model1} and \eqref{eq:model2} are generally nonconvex and thus may have multiple minimizers (of various quality),
for best results in practice one should initialize GRANSO from multiple starting points and take the best
of the resulting computed solutions.
Finally, for all problems and Approaches~1 and 2, we set
$\texttt{tol}_\mathrm{sa} = 0.9 \cdot \min \{ \alpha_\mathrm{MCK}(v_\mathrm{init}),  \alpha_\mathrm{MCK}(v_\mathrm{zero}) \}$,
where \mbox{$v_\mathrm{zero} = \texttt{zeros(3,1)}$},

\begin{figure}[t]
\centering
\subfloat[Spectral abscissa minimization only]{
\resizebox*{\figsizes}{!}{\includegraphics[trim=0cm 0cm 0cm 0cm,clip]{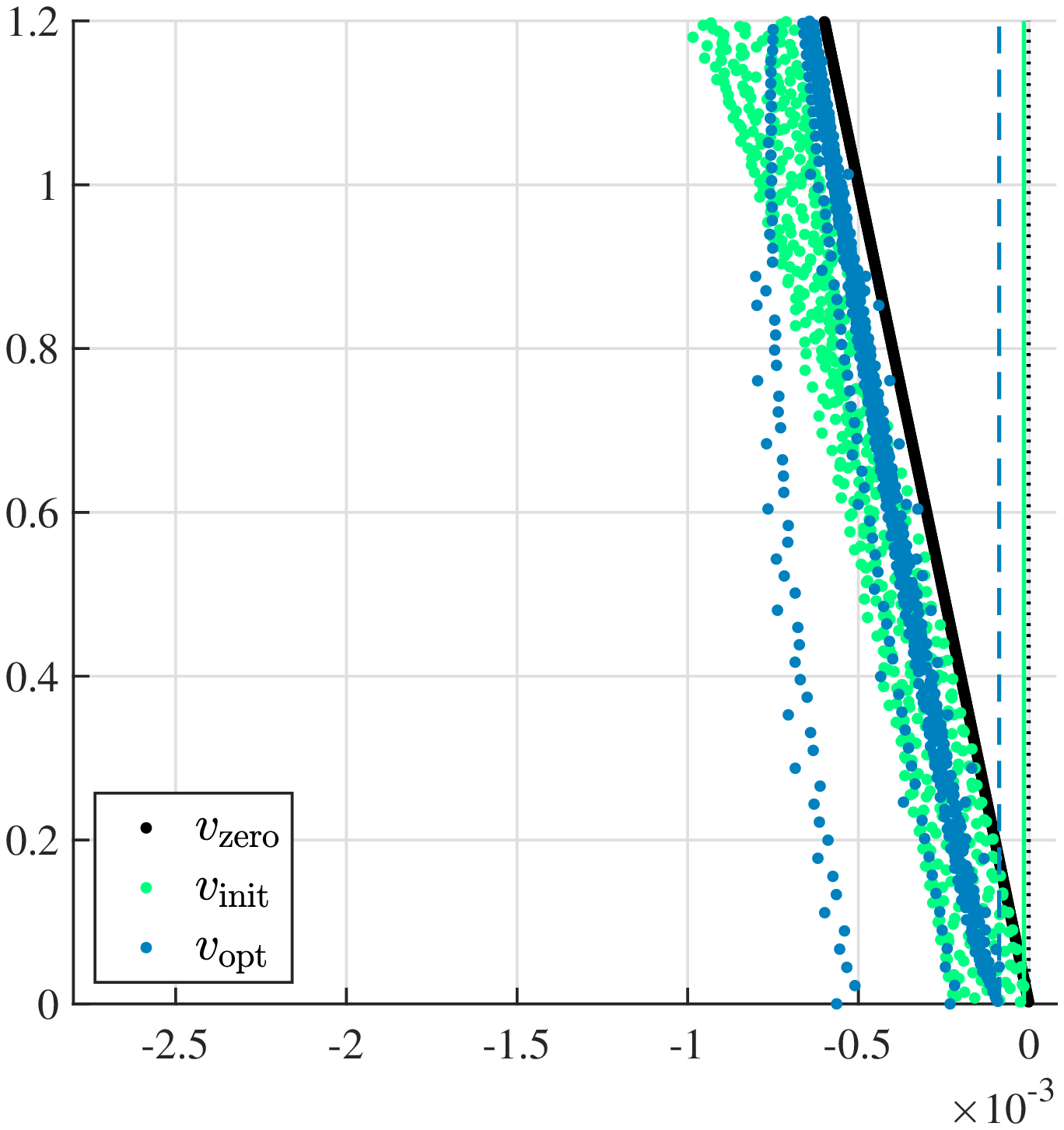}}
}
\subfloat[Approach~1 (Fixed ellipses)]{
\resizebox*{\figsizes}{!}{\includegraphics[trim=0cm 0cm 0cm 0cm,clip]{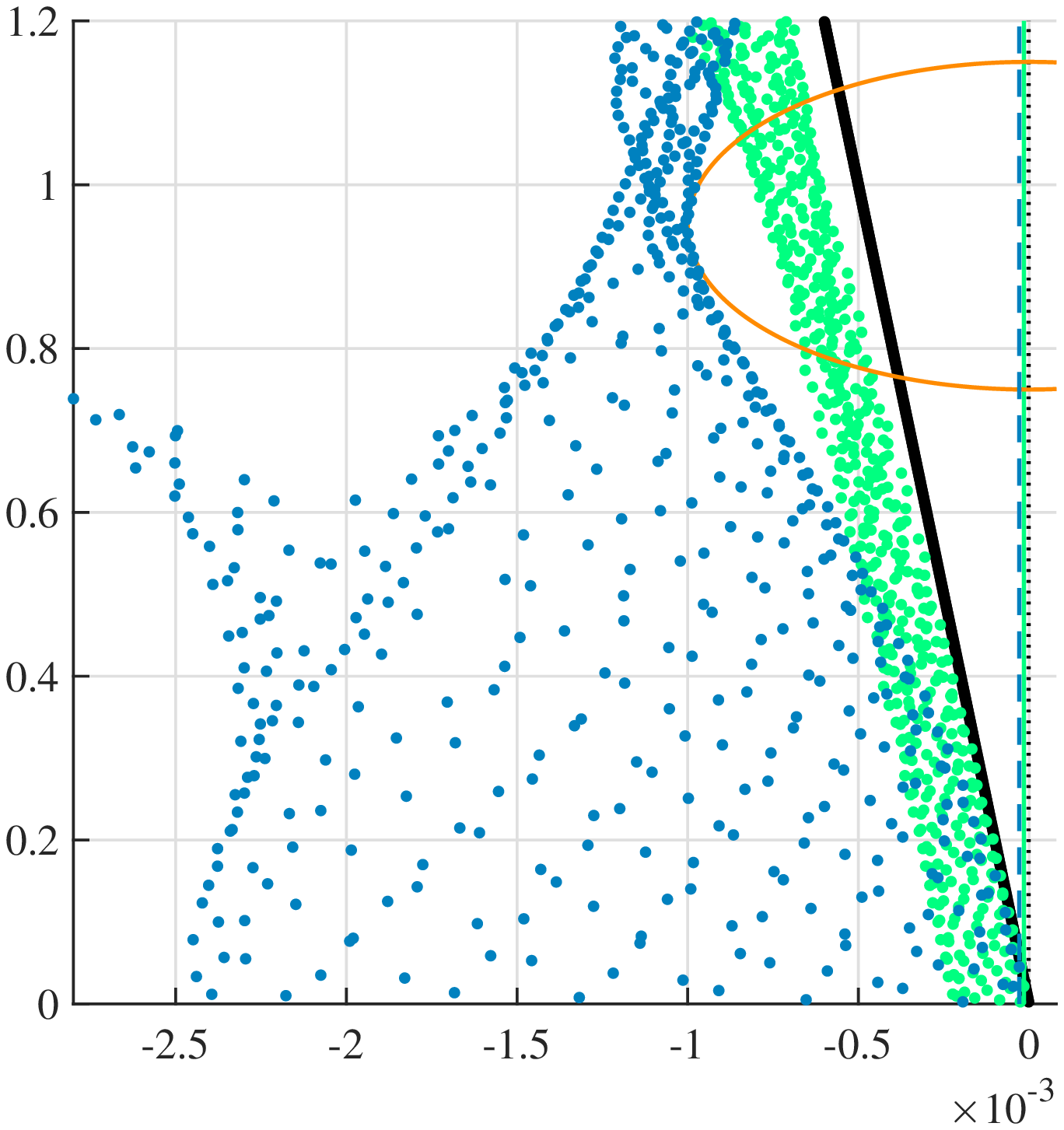}}
}
\caption{
The eigenvalues (depicted as dots) are shown for $v_\mathrm{zero}$ (no external damping),
$v_\mathrm{init}$ (the damping at the initial point), and $v_\mathrm{opt}$ (the optimized viscosities computed by GRANSO).
The spectral abscissa for each set of viscosities is depicted via a vertical line: $v_\mathrm{zero}$ (dotted),
$v_\mathrm{init}$ (solid), and $v_\mathrm{opt}$ (dashed).
For Approach~1 (right), the ellipse defining our frequency-weighting constraint $d_{\Lambda,\mathcal{E}}(v) \geq 1$,
which requires that none of the eigenvalues are inside
the ellipse defined by $E = (0.001,0.2, 0.95\imagunit)$, is also shown.
}
\label{fig:approach1}
\end{figure}

We begin with Approach~1, where we used $\alpha=0.001$ for $C_\mathrm{int}$ and defined
constraint $d_{\Lambda,\mathcal{E}}(v) \geq 1$ from \eqref{eq:model1} using a single ellipse,
specifically $E=(0.001,0.2,0.95\imagunit)$.
For comparison purposes, we also ran~\cref{alg:opt} a second time on this same problem
but without our ellipse constraint.
In~\cref{fig:approach1}, we show the different eigenvalue configurations before and after optimization.
When only minimizing the spectral abscissa, GRANSO ran for 7 iterations,
while for Approach~1, GRANSO ran for 16 iterations.  The solutions returned by GRANSO were, respectively,
\[
	v_\mathrm{opt} = \left[ \begin{array}{r} 238.7 \\ 101.2 \\ 132.6 \end{array} \right]
	\qquad \text{and} \qquad
	v_\mathrm{opt} =  \left[ \begin{array}{l} 8.117 \\ 8.187 \\ 0.001467 \end{array} \right].
\]
From \cref{fig:approach1}, we clearly see that both of these solutions are close to the nonsmooth manifold,
with the former resulting in several eigenvalues being close to attaining the spectral abscissa (the left pane)
and the latter resulting in many more eigenvalues being exceptionally close to the boundary of our specified ellipse (the right pane).
Moreover, we see that while the addition of constraint $d_{\Lambda,\mathcal{E}}(v) \geq 1$ in Approach~1 causes
the spectral abscissa to be minimized less,
Approach~1 did in fact move all of the eigenvalues at the initial viscosities $v_\mathrm{init}$ out
of our ellipse region.  In other words, Approach~1 successfully computed a feasible set of viscosities
that both selectively and significantly damped the desired frequency band.

\begin{figure}[t]
\centering
\subfloat[$\alpha = 0.004$]{
\resizebox*{\figsizes}{!}{\includegraphics[trim=0cm 0cm 0cm 0cm,clip]{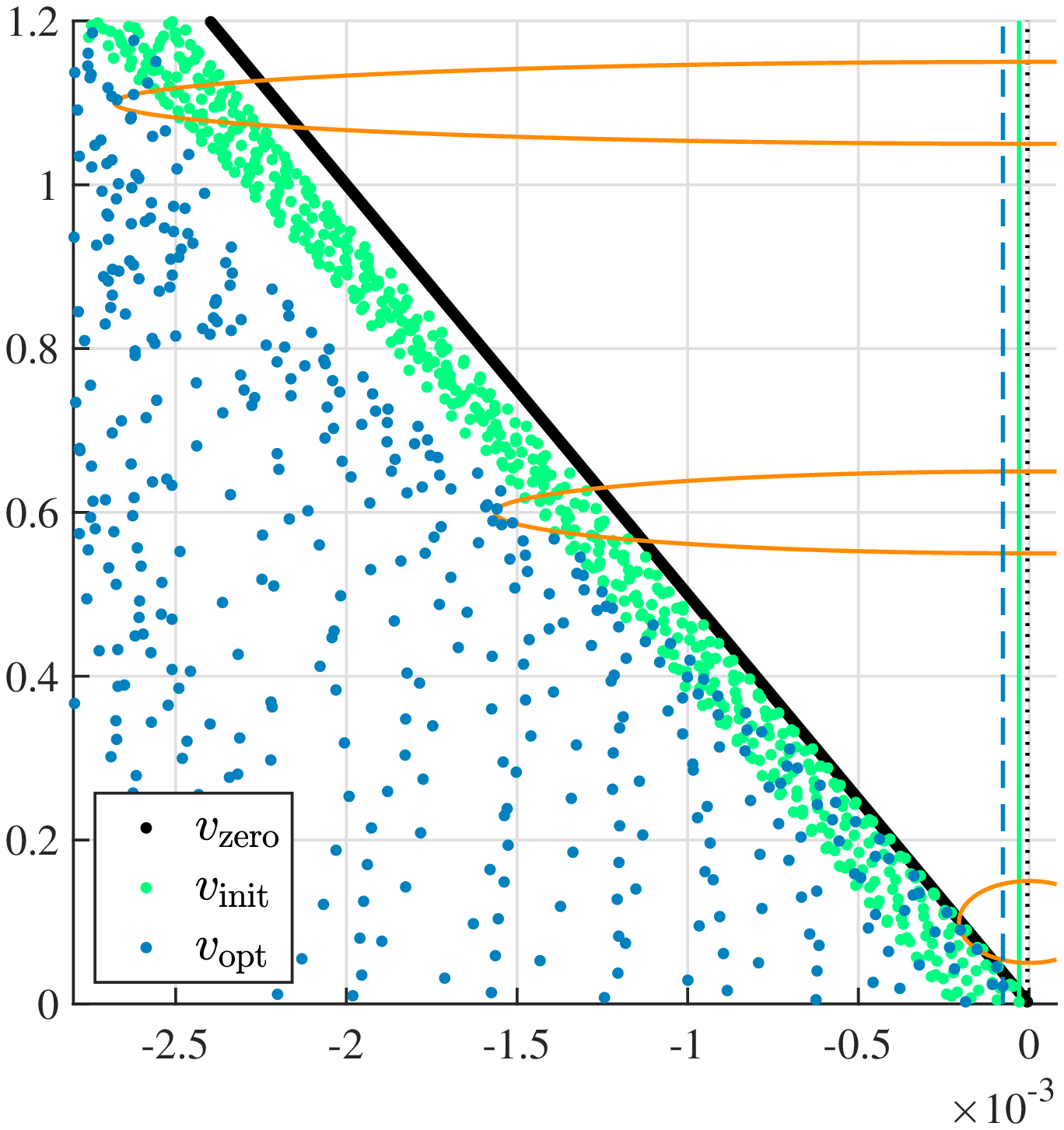}}
}
\subfloat[$\alpha = 0.0004$]{
\resizebox*{\figsizes}{!}{\includegraphics[trim=0cm 0cm 0cm 0cm,clip]{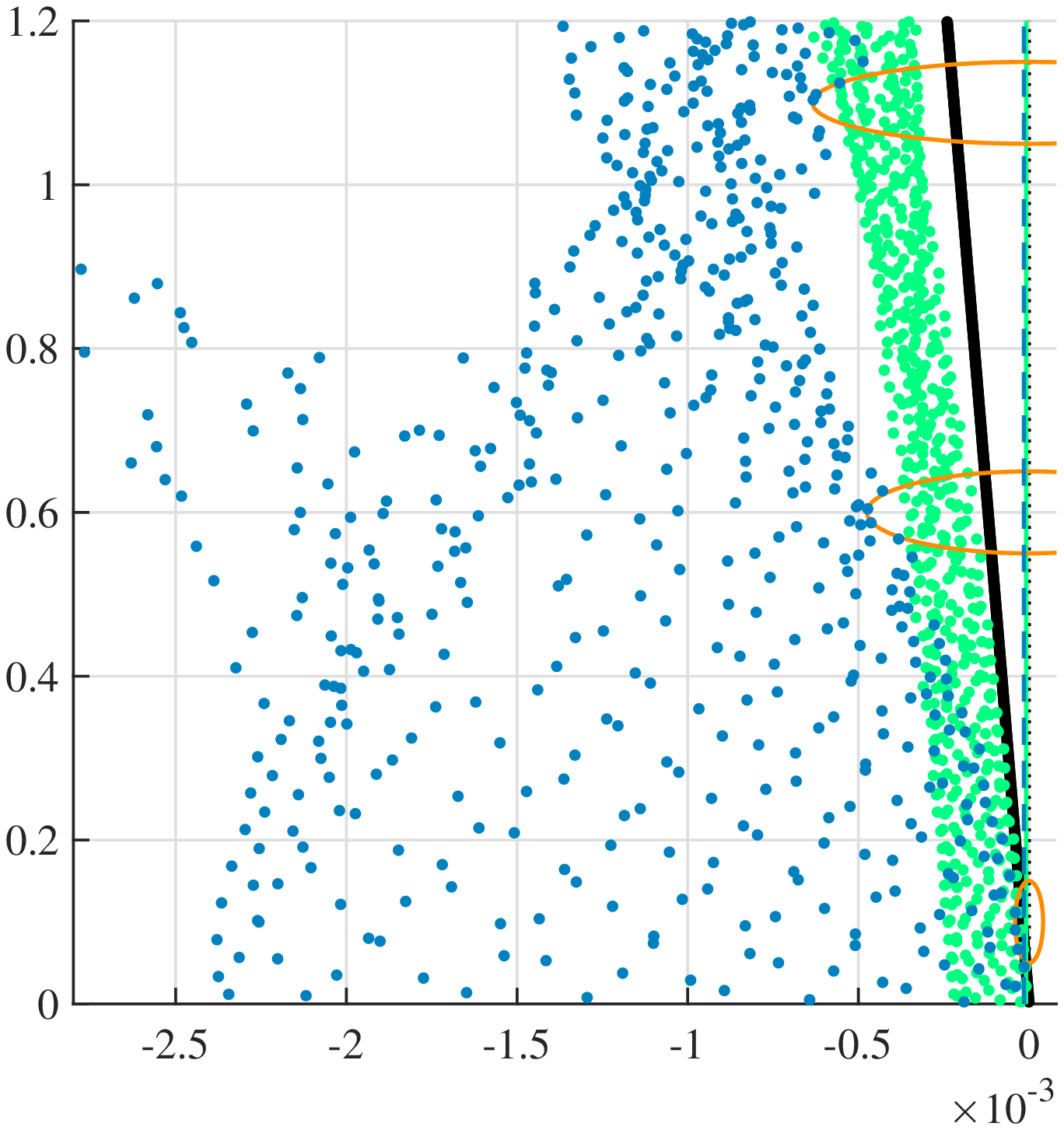}}
}
\caption{The resulting spectrum configurations and maximized ellipses computed by Approach~2 (Variable ellipses) are shown for two different choices
of $\alpha$ determining the internal damping matrix $C_\mathrm{int}$.
See the caption of~\cref{fig:approach1} for the description of the legend.
}
\label{fig:approach2}
\end{figure}

Turning to Approach~2, we used three ellipses to define our objective function in \eqref{eq:model2},
specifically $E_1 = (\sim,0.05,0.1\imagunit)$, $E_2 = (\sim,0.05,0.6\imagunit)$, and $E_3 = (\sim,0.05,1.1\imagunit)$,
with respective weightings $\phi_1 = 1$, $\phi_2 = 0.2$, and $\phi_3 = 0.1$, and $\eta = 0$.
We then ran \cref{alg:opt} using this instance of Approach~2 on the same $n$-mass oscillator example that we used to test Approach~1,
except that now we used $\alpha=0.004$ and $\alpha=0.0004$.  The configurations of eigenvalues
before and after optimization are shown in \cref{fig:approach2}.  For $\alpha=0.004$ and $\alpha=0.0004$,
GRANSO respectively ran for 32 and 27 iterations before halting and respectively returned
\[
	v_\mathrm{opt} = \left[ \begin{array}{r} 8.138 \\ 7.147 \\ 1.789 \end{array} \right]
	\qquad \text{and} \qquad
	v_\mathrm{opt} =  \left[ \begin{array}{r} 8.295 \\ 7.767 \\ 1.673 \end{array} \right].
\]
for the optimized viscosity values.
We again see that the solutions returned by GRANSO are very close to the nonsmooth manifold.
In the left pane of \cref{fig:approach2}, we see that each ellipse is essentially touching at least two eigenvalues,
while in the right pane, the three ellipses are very close to touching three, four, and three eigenvalues, respectively, from top to bottom.
Furthermore, the resulting eigenvalue configurations in \cref{fig:approach2} confirm that
Approach~2 is indeed able to perform the desired frequency-weighted damping, as specified by
the semi-minor axis values and centers of ellipses~$E_1$, $E_2$, and $E_3$.

\section*{Acknowledgment} This work has been fully supported by Croatian Science Foundation under the project `Vibration Reduction in Mechanical Systems' (IP-2019-04-6774).

\footnotesize
\bibliographystyle{plain}
\bibliography{csc,mor}

\end{document}